\documentclass[12pt]{amsart}
\usepackage{amsmath}\usepackage{amsrefs}\usepackage{ifsym}\usepackage{amssymb}\usepackage{upgreek}
\makeatletter
\newcommand{\leqnomode}{\tagsleft@true}
\newcommand{\reqnomode}{\tagsleft@false}
\makeatother
\usepackage{amsfonts}\usepackage{verbatim}
\usepackage{hyperref}
\usepackage{amssymb}\usepackage{marvosym}
\usepackage{dictsym}\usepackage{wasysym}
\usepackage{amstext}
\usepackage{amsbsy}\usepackage{bbold}
\usepackage{amsopn}\usepackage{MnSymbol}
\usepackage{amsthm}\usepackage{color}\usepackage{tikz}\usetikzlibrary{arrows}
\RequirePackage{filecontents}

\DeclareMathAlphabet{\mathpzc}{OT1}{pzc}{m}{ithttp://tex.stackexchange.com/questions/180261/math-in-table-of-contents-bookmarks-and-heading}

\def\proclaim#1{\vskip0.5em\noindent{\bf #1}\it }
\def\endproclaim{\vskip0.5em\par\noindent\rm}

\def\proclaim#1{\vskip0.5em\noindent{\bf #1}\it}
\def\endproclaim{\vskip0.5em\par\noindent\rm}

\def\undersetbrace#1#2{\underbrace{#2}_{#1}}

\def\demo#1{\vskip0.5em\noindent{\bf #1\ }}

\def\text#1{\mbox{#1}}
\def\flushpar{\par\noindent}

\def\mod{\mbox{ mod }}

\newcommand{\mapright}[1]{%
    \smash{\mathop{%
        \hbox to 1cm{\rightarrowfill}
        }
    \limits^{#1}
    }
}
\newcommand{\mapleft}[1]{%
    \smash{\mathop{%
        \hbox to 1cm{\rightarrowfill}
        }
    \limits_{#1}
    }
}

\def\e{\epsilon}
\def\a{\alpha}
\def\b{\beta}
\def\G{\Gamma}
\def\g{\gamma}
\def\d{\delta}
\def\D{\Delta}
\def\s{\sigma}
\def\Si{\Sigma}
\def\th{\theta}
\def\l{\lambda}
\def\x{\times}
\def \N{\Bbb N}\def\Z{\Bbb Z}\def\R{\Bbb R}
\def \F{\mathcal F}

\def\f{\flushpar}
\def\u{\underline}
\def\v{\varphi}

\def\om{\omega}
\def\Om{\Omega}
\def\B{\mathcal B}
\def\({\biggl(}
\def\){\biggr)}
\def\<{\langle}
\def\>{\rangle}
\def\bdy{\partial}
\def\bul{\smallskip\f$\bullet\ \ \ $}\def\lfl{\lfloor}\def\rfl{\rfloor}
\def\lcl{\lceil}\def\rcl{\rceil}
\def\){\biggr)}

\def\st{\text{such that}}

\def\<{\bold\langle}
\def\>{\bold\rangle}

\def\bul{\smallskip\f$\bullet\ \ \ $}\def\sms{\smallskip\f}\def\Lra{\Longrightarrow}\def\lra{\longrightarrow}\def\Par{\smallskip\f\P}
\def\Smi{\smallskip\f{\LARGE \smiley}}

\def\bdy{\partial}\def\vvvert{\vert\vert\vert}
\def\dist{\text{\tt dist}\,}
\begin{document}
\title{{Rational ergodicity  of 
Step function Skew Products.}}
\author{ Jon. Aaronson , Michael Bromberg $\&$ Nishant  Chandgotia}
 \address[Aaronson]{\ \ School of Math. Sciences, Tel Aviv University,69978 Tel Aviv, Israel.\f\ \ \  {\it Webpage }: {\tt http://www.math.tau.ac.il/$\sim$aaro}}\email{aaro@post.tau.ac.il}
\address[Bromberg]{\ \ Math. Dept., Bristol University, Bristol, UK.}\email{micbromberg@gmail.com}

\address[Chandgotia]{ School of Math. Sciences, Tel Aviv University,69978 Tel Aviv, Israel.}\email[Chandgotia]{nishant.chandgotia@gmail.com}
\thanks { The research of Aaronson and Chandgotia was partially supported by ISF grant No. 1599/13.   
Chandgotia was also partially supported by ERC grant No. 678520 $\&$  Bromberg's research was  supported by
 EPSRC Grant EP/I019030/1.\ \ \copyright  2016-7. }

\begin{abstract}We study {\tt rational step function} skew products over certain rotations of the circle proving ergodicity and bounded rational ergodicity 
when  rotation number is a quadratic irrational. The latter arises from a consideration of the asymptotic temporal statistics of an 
orbit  as modelled by an {\tt associated affine random walk}.
\end{abstract}
\maketitle\markboth{Step function Skew Products.}{Aaronson, Bromberg,  $\&$ Chandgotia}

   \section*{\S0 Introduction}

A {\it rational step function} is a right continuous, step function on the additive circle $\Bbb T:=\Bbb R/\Bbb Z\cong [0,1)$ taking values in $\Bbb R^d$, whose 
discontinuity points are rational. 

\

Let $\v:\Bbb T\to\Bbb R^d$ be a rational step function.

\

The {\it skew products} $T_{\a,\v}=T_\a:\Bbb T\x\Bbb R^d\to\Bbb T\x\Bbb R^d\ \ (\a\in\Bbb T)$ defined by
$$T_{\a,\v}(x,y):=(x+\a,y+\v(x))$$

are conservative if and only if
\begin{align*} \int_\Bbb T\v(t)dt=0.\end{align*}

Necessity follows from the ergodic theorem and sufficiency follows from the Denjoy-Koksma inequality (see below).

Consider the collections of {\tt badly approximable} irrationals
$$\text{\tt BAD}:=\{\a\in\mathbb R\setminus\mathbb Q:\ \exists\ \th>0,\ |\a-\tfrac{p}q|\ge\tfrac{\th}{q^2}\}$$ and of quadratic irrationals
$$\text{\tt QUAD}:=\{\a\in\mathbb R\setminus\mathbb Q:\ \a \ \ \text{quadratic}\}.$$
It is known  that $\text{\tt QUAD}\subset\text{\tt BAD}$ and that $\text{\tt BAD}$ has Lebesgue measure zero (see e.g. \cite{HW}).

\

\subsection*{Denominator of a rational step function}
\

Fix $d,\ Q\in\Bbb N,\ Q\ge 2$ and $\Phi:\Bbb Z_Q\to\Bbb R^d$.

The {\it rational step function with denominator $Q$ and values $\Phi$} is the step function $\v=\v^{(\Phi)}:\Bbb T\to\Bbb R^d$ defined by
$$\v(x)=\Phi(\upkappa(x)),$$ where $\upkappa:[0,1)\to\Bbb Z_Q$ is defined by $\upkappa(x):=\lfl Qx\rfl$.  Every rational step function is of this form for some $Q\ge 2$ and $\Phi:\Bbb Z_Q\to\Bbb R^d$.

If $\v:\Bbb T\to\Bbb R^d$ is a rational step function with denominator $Q$, then
$$\int_\Bbb T\v(x)dx=\frac1Q\sum_{k=0}^{Q-1}\Phi(k).$$
We prove, for $\v:\Bbb T\to\Bbb R^d$ a rational step function with denominator $Q$ and values $\Phi$ and which is {\it centered} in the sense that $\int_\Bbb T\v(x)dx=0$:
\proclaim{Theorem 1':\ \ Ergodicity}
\

There is a collection $\text{\tt SBAD}\subset\Bbb R\setminus\Bbb Q$ of full 
Lebesgue measure so that $\text{\tt SBAD}\supset\text{\tt BAD}$ and so that
if $\a\in\text{\tt SBAD}$, then $(\Bbb T\x\G ,\B(\Bbb T\x\G ),m_\Bbb T\x m_\G ,T_{\a,\v})$ is a {\tt CEMPT} where
$\G :=\overline{\<\v(\Bbb T)\>}$ is the closed subgroup of $\Bbb R^d$ generated by $\Phi(\Bbb Z_Q)$.\endproclaim
\

Here and throughout, {\tt CEMPT} means  conservative, ergodic, measure preserving transformation, $m_{\Bbb G}$ denotes Haar measure on the locally compact, Polish, Abelian group $\Bbb G$, normalized if $\Bbb G$ is compact. Also, $\v:\Bbb T\to\Bbb R^d$ is always going to mean a centered rational step function.
\

Theorem 1' will follow from the stronger theorem 1 (see below).
\

The technique  of the proof of theorem 1 is not new. 
For older, related results, see \cite{CP},\ \cite{Oren} and references therein.

\proclaim{Theorem 2:\ \ Temporal CLT}\ 

\

\ If  $\a\in\text{\tt QUAD}$ 
and $\dim\text{\tt span}_\Bbb R\,\v(\Bbb T)=d$, then $\exists\ \ell_k\in\Bbb N,\ \ell_k\uparrow\ \&\ \ell_k\propto\l^k$ 
for some $\l>1$ and $\mu^{(0)}\in\Bbb R^d$
so that for any box $I\subset\Bbb R^d$,
$$\frac1{\ell_k}\#\,\{1\le n\le \ell_k:\ \frac{\v_n(0)-k\mu^{(0)}}{\sqrt{k}}\in I\}\xrightarrow[k\to\infty]{}\int_If_Z(t)dt$$ where $Z$ is a 
globally supported, centered, normal random variable on $\Bbb R^d$ and $f_Z$ is its probability density function.\endproclaim
Here and throughout $\v_n(x):=\sum_{k=0}^{n-1}\v(x+k\a)$ and $\#$ denotes counting measure.
\

For an introduction to
temporal statistics in dynamics see \cite{DS}. 
Theorem 2 here  is a generalization of a subsequence version of theorem 1.1 in \cite{Beck}, which in turn has been recently strengthened in \cite{BU}.
\proclaim{Theorem 3:\ \ Rational ergodicity}\ \ \ \ 

\

Suppose that\  $\a\in\text{\tt QUAD}$\  and that  $\<\v(\Bbb T)\>=\Bbb Z^d$, 
then $$(\Bbb T\x\Bbb Z^d,\B(\Bbb T\x\Bbb Z^d),m_\Bbb T\x\#,T_{\a,\v})$$ is
{ boundedly rationally ergodic} and $a_n(T_{\a,\v})\asymp\frac{n}{(\log n)^\frac{d}2}$.\endproclaim
See \cite{AK} for a definition of {\tt bounded rational ergodicity}. Bounded rational ergodicity of $T_{\a,\v}$ for $\v=1_{[0,\frac12)}-1_{[\frac12,1)}$
was established in \cite{AK} for $\a\in\text{\tt QUAD}$ and in \cite{ABN2016} for $\a\in\text{\tt BAD}$.
\

\subsection*{Notations} Here and throughout, for $a_n,\ b_n>0$:
\sms $a_n\ll b_n$ means $\exists\ M>0$ so that $a_n\le Mb_n$ for each $n\ge 1$,
\sms $a_n\asymp b_n$ means $a_n\ll b_n$  and $b_n\ll a_n$,
\sms $a_n\propto b_n$ means $\exists\ \lim_{n\to\infty}\frac{a_n}{b_n}\in \Bbb R_+:=(0,\infty)$ and
\sms $a_n\sim b_n$ means $\frac{a_n}{b_n}\xrightarrow[n\to\infty]{}\ 1$.

\subsection*{Outline of the rest of the paper. }
\

In \S1 we prove theorem 1, a stronger version of theorem 1'.
 The rest of the paper is devoted to the proofs
of theorems 2 and 3.
\

As in \cite{AK} , \cite{ABN2016},  proofs
rely on recursive properties of the tuples $\left(\varphi_{n}\left(0\right):1\leq n\leq \ell_{k}\right)$,
for a suitably chosen sequence $\ell_{k}\uparrow\infty$. 
\

To study the temporal statistics of these tuples, we consider the ``temporal random variables'' $x_{k}:\left\{ 1,...,\ell_{k}\right\} \rightarrow\mathbb{R}^{d}$,
defined by $x_{k}\left(n\right)=\varphi_{n}\left(0\right)$, where
 $n$ is a uniformly
distributed random variable with values in $\left\{ 1,...,\ell_{k}\right\} $. In other words, 
$$\text{\tt Prob}\,\left(x_{k}\in I\right)=\frac{1}{\ell_{k}}\#\left\{ 1\leq n\leq \ell_{k}:\ \varphi_{n}\left(0\right)\in I\right\}.$$

The recursive properties of the tuples (see \S2),
 allow us to construct an associated affine random walk ({\tt ARW}) 
 which models the distribution of the ``temporal random variables'' (see \S3).

\

In \S4 we show    that when $\alpha$ is quadratic,
 the sequence of expectations $E\left(x_{k}\right)$
is asymptotically linear. This culminates in the approximation of the
distribution of $x_{k}-E\left(x_{k}\right)$ by an affine random walk
generated by a sequence of centered, independent, identically distributed affine transformations
(see the {\tt ARW} centering lemma). 
\

This enables proof in \S5 of theorem 2 which is a central limit theorem for $(x_k:\ k\ge 1)$.
The proof of theorem  3
in  \S6 is based on a ``weak, rough local limit theorem'' for $(x_k:\ k\ge 1)$. Both proofs  use a spectral theory of 
{\tt ARW}s based on perturbation theory of stochastic matrices (as in \cite{HH}).
\section*{\S1 Ergodicity}
\

\subsection*{ Regular continued fractions}
\

 Recall that the {\it regular continued fraction expansion} of $\a\in (0,1)\setminus\Bbb Q$ is 
\begin{align*}\a & =\frac{1}{a_1+\frac{1}{a_2+_{\ddots+\frac{1}{a_{n}+_{\ddots}}}}}\\ & \,
\\ &=:\cfrac[r]{1|}{|a_1}+\cfrac[r]{1|}{|a_2}+\dots+\cfrac[r]{1|}{|a_{n}} \ \ +\cdots\\ &=
(a_1,a_2,\dots)\end{align*}
where $a_n:=a(G^{n-1}\a)\in\Bbb N$ with $a(\a ):=\lfl\tfrac1\a \rfl\ \&\ 
 G(\a ):=\{\tfrac1\a \}=\alpha-\lfloor \frac{1}{\alpha}\rfloor$ for $\a \in\Bbb T\setminus\Bbb Q$.

\

Recall  that $G((0,1)\setminus\Bbb Q)\subset(0,1)\setminus\Bbb Q$ and so  every irrational in $(0,1)$ indeed has an infinite regular continued fraction expansion.
On the other hand, if $\a\in (0,1)\cap\Bbb Q$ then $\exists\ n\ge 1,\ G^n(\a)=0$ and $\a$ has only a finite regular continued fraction expansion. 
In the sequel, we'll consider modified continued fractions
where the situation is different.
\par Fix $\a=(a_1,a_2,\dots)\in (0,1)\setminus\Bbb Q$ and $n\ge 1$ and define the  {\it principal convergents} 
$\frac{p_n}{q_n},\ p_n,\ q_n\in\Bbb N_0,\ \gcd\,(p_n,q_n)=1$ by
$$\frac{p_n}{q_n}:=\cfrac[r]{1|}{|a_1}+\cfrac[r]{1|}{| a_2}+\dots+\cfrac[r]{1|}{|a_n }.$$
 Here, and throughout, for $k\in\Bbb Z,\ k\ge 0$,  we denote
 $$\Bbb N_k:=\{n\in \Bbb Z:\ n\ge k\}.$$
 
 \
 
 The {\it principal denominators} $q_n$ of $\a $ are given by
$$q_{0}=1,\ q_{1}=a_1,\ q_{n+1}=a_{n+1}q_{n}+q_{n-1};$$
the numerators $p_n$  are given by
$$p_{0}=0,\ p_{1}=1,\ p_{n+1}=a_{n+1}p_{n}+p_{n-1}$$
and the {principal convergents} $\tfrac{p_n}{q_n}$ satisfy
$$\frac1{q_n(q_n+q_{n+1})}<|\a-\frac{p_n}{q_n}|<\frac1{q_nq_{n+1}}.$$
We'll also need theorems 16, 17 and 19 in \cite{Khintchine}:
\proclaim{Proposition}\ \ For $\a\in\Bbb T\setminus\Bbb Q$,
\sms If for some $0\leq a\leq b (b\in \N)$, $|b\a-a|<|d\a-c|$ for all $0<d<b$ then $b=q_k$ for some $k\in \N$.
\sms For $0\le a\le b<q_k\ \ (b,k\in \N)$, $|b\a-a|>|q_k\a-p_k|$.
\sms If $p,q\in\Bbb N,\ (p,q)=1$ and $|\a-\tfrac{p}q|<\frac1{2q^2}$, then $q\in\{q_k:\ k\ge 1\}$.\endproclaim

The following is also well known (see e.g. \cite{HW},\ \cite{Khintchine}):
 
\proclaim{Proposition} \ \ \ Let $\a=(a_1,a_2,\dots)\in (0,1)\setminus\Bbb Q$, then
\sms {\rm (i)} $\a\in\text{\tt QUAD}$ iff $\exists\ K,\ L\ge 1$ so that $a_{k+L}=a_k\ \forall\ k\ge K$;
\sms {\rm (ii)} $\a\in\text{\tt BAD}$ iff  $\sup_{k\ge 1}\,a_k<\infty$.\endproclaim
\

For $Q\ge 2$, we'll also need the  collection
 $$\text{\tt SBAD}_Q:=\{\a:\  \varlimsup_{n\to\infty} \tfrac{q_{n'}}{q_{n+1}}>0\ \text{ where }n':=\max\{1\le m<n:\ \tfrac{q_m}{q_n}<\tfrac1Q\}\}.$$
Evidently, $\text{\tt BAD}\subset \text{\tt SBAD}:=\bigcap_{Q\ge 2}\text{\tt SBAD}_Q$ and it is not hard to show that $\text{\tt SBAD}$   has full Lebesgue measure.
\

In the following, $\v=\v^{(\Phi)}:\Bbb T\to\Bbb R^d$ is a rational step function with denominator $Q\ge 2$.
\proclaim{Theorem 1}
\

Suppose that either {\rm (i)}  $\a\in\text{\tt SBAD}_Q$, or {\rm (ii)} $\a\notin\Bbb Q\ \&\ Q$ is prime, then
$(\Bbb T\x\G ,\B(\Bbb T\x\G ),m_\Bbb T\x m_\G ,T_{\a,\v})$ is a {\tt CEMPT}.\endproclaim
The rest of this section is devoted to the proof of theorem 1.
\subsection*{
  Essential values and Periods}

 Let $(X,\B,m)$ be a standard probability space,
 and let $T:X\to X$ be an invertible, ergodic, probability preserving transformation and $\B_+:=\{A\in \B:m(A)>0\}$.
 \

Suppose that   $\Bbb G$ is a locally compact, Polish, Abelian group equipped with the translation invariant metric  $\rho$ (e.g. $\rho(x,y)=\|x-y\|$ if $\Bbb G\le\Bbb R^d$).

Let  $\v:X\to \Bbb G$ be
measurable and define $\v _n: X\longrightarrow\Bbb G$ by 
$$\v _n:=\sum_{k=0}^{n-1}\v\circ T^k.$$

\

The  collection of {\it  essential values}    of $\v$ (as in \cite{rigveda}) is 
$$E(\v):=\{a\in \Bbb G:\forall\ A\in\B_+,\ \e>0,\ \exists\ n\in\Bbb Z,\
m(A\cap T^{-n}A\cap [\rho(\v_n,a)<\e])>0\}.$$
\

The {\it skew product} $T_\v:X\x\Bbb G\to X\x\Bbb G$ is defined by
$$T_\v(x,y):=(Tx,y+\v(x))$$
and preserves the measure $m\x m_{\Bbb G}$.

\

\par Define the collection of {\it periods} for $T_\v$-invariant functions:
$$\text{\tt Per}\,(\v)=\{a\in \Bbb G: \tau_a A=A\ \mod\ m\times m_{\Bbb G}\ \forall\ A\in\B(X\x\Bbb G),\  T_\v(A)=A\}$$
where $ \tau_a(x,y)=(x,y+a).$
\

It is not hard to see that $T_\v$ is ergodic iff $T$ is ergodic $\&\ \text{\tt Per}\,(\v)=\Bbb G$.
\proclaim{Schmidt's Theorem\ \ \cite{rigveda}}
\par  
 $E(\v)$ is a closed subgroup of $\Bbb G$ and $E(\v)=\text{\tt Per}\,(\v)$.
\endproclaim

In view of this, the conclusion of  theorem 1 is equivalent to
\begin{align*}\tag{\Wheelchair}\label{wheelchair}\G:=\overline{\Phi(\Bbb Z_Q)}=E(\v). 
\end{align*}

\

We prove this first in the case that $\overline{\Phi(\Bbb Z_Q)}$ is countable and then deduce the uncountable case. 
\

Let

$$D(\a):=\{q\in\Bbb N:\ \exists\ p\in\Bbb N,\ |\alpha-\tfrac{p}{q}|<\tfrac{1}{q^2}\}.$$
We'll need

\proclaim{ Denjoy-Koksma inequality \ \ \ (\cite{Katznelson},\ \cite{Herman})}\ \ 
$$\|\v_q\|_\infty\le\bigvee_{\Bbb T}\v\ \forall\ q\in\ D(\a),$$
where $\bigvee_{\Bbb T}\v$ denotes the total variation of $\v$.
\endproclaim

\subsection*{Remark}
Consequently, when $\G=\overline{\<\Phi(\Bbb Z_Q)\>}$ is countable, there is a finite set $F\subset \G$ such that $\v _q(x)\in F$ for every $x\in\Bbb T\ \&\ q\in D(\a)$.

Given an Abelian group $\Bbb G$ and $g_0\in \Bbb G$, let $r_{g_0}: \Bbb G\longrightarrow \Bbb G$ denote the group rotation on $\Bbb G$ given by $r_{g_0}(g):=g+g_0$.
\demo{Proof of theorem 1 in the countable case}
\proclaim{Sublemma 1}  \ \ For theorem 1 in the countable case, it suffices that 
\begin{align*}
\tag{\Gentsroom}\label{Gentsroom}\Phi(\epsilon+1)-\Phi(\epsilon)\in \text{\tt Per}\,(\v)\ \forall\ \epsilon\in \Bbb Z_Q.
\end{align*}
\endproclaim\demo{Proof}
\

Let $\G_0\subset \G$ be the group generated by $\{\Phi(\epsilon+1)-\Phi(\epsilon)~:~\epsilon\in \Bbb Z_Q\}$, that is, 
$$\G_0:=\<\{\Phi(\epsilon+1)-\Phi(\epsilon)~:~\epsilon\in \Bbb Z_Q\}\>\ \le\ \G.$$
Evidently, $\Phi(\epsilon)+\G_0=\Phi(0)+\G_0\ \forall\ \epsilon\in\Bbb Z_Q$ whence
$\v+\G_0\equiv\Phi(0)+\G_0$ and $\G/\G_0$ is cyclic.
\

We claim moreover that $\#\G/\G_0\le Q$. To see this, 
using $\sum_{\epsilon\in\Bbb Z_Q}\Phi(\epsilon)=0$, we have 
$$\G_0\ni\sum_{\epsilon\in\Bbb Z_Q}(\Phi(0)-\Phi(\epsilon))=Q\Phi(0)$$
whence indeed  $\#\G/\G_0\le Q$.
\

By (\Gentsroom), $\G_0\subset \text{\tt Per}\,(\v)$. 
\

By Schmidt's theorem, if $h\in L^\infty(\Bbb T\x\G)$ and $h\circ T_{\a,\v}=h$, then
$h\circ\tau_a=h$ a.e. $\forall\ a\in\G_0$ and $\exists\ H\in L^\infty(\Bbb T\x\G/\G_0)$ so that
$$h(x,\g)=H(x,\g+\G_0)\ \text{for a.e.}\ (x,\g)\in \Bbb T\x\G.$$
Evidently
$$H\circ T_{\a,\psi}=H\ \text{ a.e.}$$
where $\psi:\Bbb T\to\G/\G_0,\ \psi:=\v+\G_0\equiv\Phi(0)+\G_0$ (as before).
\

Defining  $T_{\a,\psi}:\Bbb T\x\G/\G_0\to \Bbb T\x\G/\G_0$  as usual, we have
\

$$T_{\a,\psi}\cong r_\a\x r_{\Phi(0)+\G_0}:\Bbb T\x\G/\G_0\to \Bbb T\x\G/\G_0.$$
which
is  ergodic, being a product of two ergodic group rotations with disjoint spectra.

\

Thus $H$ is constant a.e., whence also $h$, and $T_{\a,\v}$ is ergodic.\ \ \CheckedBox

\proclaim{Sublemma 2}  \ \ 
\begin{align*}
\tag{\Gentsroom}\Phi(\epsilon+1)-\Phi(\epsilon)\in \text{\tt Per}\,(\v)\ \forall\ \epsilon\in \Bbb Z_Q.
\end{align*}
\endproclaim\demo{Proof}
\

We'll prove the sublemma using  
\proclaim{Oren's Lemma  \cite{Oren} }
If there exist $n_k\in \Bbb N$ and $A_{k}\subset \Bbb T$ such that $\v _{n_k}$ is constant on $A_k$ and $\v _{n_k}|_{A_{k}}\lra a$, $\inf m_{\Bbb T}(A_{k})>0$ and $\lim_{k\longrightarrow \infty}\vvvert n_k \alpha\vvvert=0$ then $a\in \text{\tt Per}\,(\v )$.
\endproclaim
Here and and throughout,  $\vvvert  x\vvvert :=\min_{k\in\Bbb Z}|x-k|$.
\

Note that a  version of Oren's lemma is implicit in  \cite{Conze76}.

\

Next, we claim that for (\Gentsroom), it suffices to show
\Smi \ \ For any $\e \in \Z_Q$ there are sequences of measurable sets $(A_k), (B_k)\subset \Bbb T$ and positive integers $n_k\in \N$ such that $\v _{q_{n_k}}$ is constant on $A_k$ and $B_k$,
$$\v _{q_{n_k}}|_{B_k}-\v_{q_{n_k}}|_{A_k}=\Phi(\e +1)-\Phi(\e )$$
and $m(A_k), m(B_k)>c>0$ where $c$ does not depend on $k$.
\

Indeed, by the remark after  Denjoy-Koksma inequality,  there is a finite set $F$ so   that $\v _{q_{n_k}}(x)\in F \ \forall\ k\ge 1,\ x\in\Bbb T$.
\

Thus $\exists\ f\in F\ \&\ k_\ell\to\infty$ such that $\v _{q_{n_{k_\ell}}}|_{A_k}=f \ \forall\ \ell\ge 1$ whence 
$\v _{q_{n_{k_\ell}}}|_{B_k}=f+\Phi(\e +1)-\Phi(\e ) \ \forall\ \ell\ge 1$.
\

By Oren's lemma, 
$f, f+\Phi(\e +1)-\Phi(\e )\in \text{\tt Per}\,(\v )$, whence, since $\text{\tt Per}\,(\v )$ is a group,  $\Phi(\e +1)-\Phi(\e )\in \text{\tt Per}\,(\v )$ and sufficiency of {\large\smiley} is established.

\

Finally, we  construct the  sequences of measurable sets $A_k, B_k\subset \Bbb T$ as in {\large\smiley}.

To this end, we prove first that the discontinuities of $\v$ are ``dynamically separated''.
\

Let $q\in \N$. Since $\v$ is a step function with the set of discontinuities contained in $\{\frac{\ell}{Q}~:~0\leq \ell \leq Q-1\}$ the set of discontinuities of $\v_{q}$ is contained in the set
$$\left\{\frac{\ell}{Q}-\text{\tt j}\alpha:~0\leq \ell\leq Q-1\text{ and }0\leq {\text{\tt j}}\leq q-1\right\}\subset\Bbb T.$$

Hence the distance between the discontinuities is  bounded below by

$$\text{\tt disc}(q)=\min_{|\ell|\leq Q-1, |{\text{\tt j}}|\leq q-1,\ (\ell,\text{\tt j})\ne  (0,0)}\vvvert \tfrac{\ell}{Q}-{\text{\tt j}} \alpha\vvvert. $$
                                             
\proclaim{Claim}
\begin{align*}\tag{\Ladiesroom}\label{Ladiesroom}\exists\ n_k\uparrow\infty\ \&\ \th>0\ \text{s.t.}\ \text{\tt disc}(q_{n_k})\ge\frac{\th}{q_{n_k}}\ \ \forall\ k\ge 1. 
\end{align*}\endproclaim

\demo{Proof of (\Ladiesroom) when $\a\in\text{\tt SBAD}$}

By definition of {\tt SBAD}, there exist a sequence $(m_k)\in \Bbb N$, $\nu \in \Bbb N$ and $\varepsilon>0$ such that 
$$\varepsilon\ <\ \frac{q_{m_k-\nu}}{q_{m_k+1}}\ <\ \frac{q_{m_k-\nu}}{q_{m_k}}<\frac{1}{Q}$$ whence for all $r\in \Bbb Z$ and $|l|\leq q_{m_k-\nu}<\tfrac{q_{m_k}}{Q}$, since $q_{m_k}$ is a 
principal denominator,
\begin{equation*}
|\tfrac{r}{Q}-l \alpha|=\tfrac{1}{Q}|r-Ql \alpha|>\tfrac{1}{Q}|p_{m_k}-q_{m_k} \alpha|>\tfrac{1}{Q (q_{m_k}+q_{m_k+1})}>\tfrac{1}{2Qq_{m_k+1}}>\tfrac{\varepsilon}{2Qq_{m_k-\nu}}.
\end{equation*}

 (\Ladiesroom) follows. \Checkedbox
\demo{Proof of (\Ladiesroom) when $Q$ is prime}
This further splits into two separate cases.
\sms (i)
\emph{There are only finitely many $n$'s such that $q_n=0\mod Q$}: Choose $N$ large enough such that $q_n\neq 0\mod Q$ for $n\geq N$. For $n>N$,
$$\text{\tt disc}(q_{n})\ge\frac{1}{Q}\min_{0<{\text{\tt j}}<q_{n}}\{\vvvert  {\text{\tt j}}Q\alpha\vvvert , \vvvert  {\text{\tt j}} \alpha\vvvert \}.$$

As before, we have
\begin{eqnarray*}
\min_{0<{\text{\tt j}}<q_n}\vvvert  {\text{\tt j}} \alpha\vvvert >\frac{1}{2q_{n}}.
\end{eqnarray*}•

Since $q_h$ is prime to $Q$ for all $h\ge n$, for $0<{\text{\tt j}}<q_n$, $\text{\tt j}Q$ is not a multiple of $q_{n+r}$ for $r\geq 0$. Thus by \cite[theorem 19]{Khintchine} if
$$\vvvert Q{\text{\tt j}} \alpha\vvvert \leq \frac{1}{2Q{\text{\tt j}}}$$ then
$Q {\text{\tt j}}$ is a multiple of $q_r$ for some $r<n$; in this case
$$\vvvert  Q{\text{\tt j}} \alpha\vvvert  \geq \frac{1}{2q_{r+1}}\geq \frac{1}{2q_{n}}.$$
Therefore
\begin{eqnarray*}
\min_{0<{\text{\tt j}}<q_n}\vvvert  Q{\text{\tt j}} \alpha\vvvert  \geq \min \left( \vvvert  q_{n-1}\alpha\vvvert  ,\frac{1}{2Qq_n}\right)=\frac{1}{2Q q_n}\text{ implying }\text{\tt disc}(q_n)\geq \frac{1}{2Q^2 q_n}.
\end{eqnarray*}
\sms (ii)\emph{There are infinitely many $n$'s such that $q_n=0\mod Q$:} Let $(n_k)$ be the subsequence such that $q_{n_k}=0\mod Q$. Let the $\nu$-th term of the continued fraction expansion of $\alpha$ be given by $a_\nu$; we know
\begin{equation*}
\left(\begin{matrix}a_n&1\\
1&0
\end{matrix}\right)\left(\begin{matrix}a_{n-1}&1\\
1&0
\end{matrix}\right)\cdots \left(\begin{matrix}a_1&1\\
1&0
\end{matrix}\right)\left(\begin{matrix}1\\ 0
\end{matrix}\right)=\left(\begin{matrix}q_{n}\\
q_{n-1}
\end{matrix}\right).\label{equation: recursion for qns}
\end{equation*}

Since $\det\left(\begin{smallmatrix}a_\nu &1\\
1&0
\end{smallmatrix}\right)=-1$ for all $\nu$, there determinant of the product is either $1$ or $-1$. Thus $q_{n_k+1}\neq 0 \mod Q$. 
\

By the recursion formula for the principal denominators,
 we have that
$q_{n_k+r+1}=M(r)q_{n_k}+ S(r)q_{n_k+1}$ for some $M(r), S(r)\in \N$.
Again,
\begin{eqnarray*}
\min_{0<{\text{\tt j}}<q_{n_k+1}}\vvvert  {\text{\tt j}} \alpha\vvvert  >\frac{1}{2q_{n_k+1}}.
\end{eqnarray*}
If
\begin{equation*}
\vvvert  {\text{\tt j}} Q\alpha\vvvert < \frac{1}{2{\text{\tt j}}Q}\text{ for some }0<{\text{\tt j}}<q_{n_{k}+1}
\end{equation*}
then by \cite[theorem 19]{Khintchine}, ${\text{\tt j}}Q$ is a multiple of $q_\nu$ for some $\nu$. 

Since $q_{n_{k}+1}\neq 0\mod Q$, if ${\text{\tt j}}Q=\text{\tt i} q_{n_k+1}$ for some $\text{\tt i}\in \N$ then ${\text{\tt j}}\geq q_{n_k+1}$; 
it follows that ${\text{\tt j}}Q$ is not a multiple of $q_{n_{k}+1}$. Therefore $\text{\tt i}\neq n_{k}+1$. Since $q_{n_k}= 0 \mod Q$ if
$$\text{\tt j}Q=\text{\tt i} q_{n_k+r+1}=\text{\tt i} M(r)q_{n_k}+ \text{\tt i} S(r)q_{n_k+1}=0\mod Q\text{ for some }\text{\tt i} $$
then $\text{\tt i} S(r)$ is multiple of $Q$ implying ${\text{\tt j}}\geq q_{n_k+1}$; thus the number ${\text{\tt j}}Q$ cannot be a multiple of $q_{n_k+r+1}$ for any $r\in \N$ and it 
follows that $\text{\tt i}\leq n_k$.

Hence  by \cite[theorem 16]{Khintchine} we have
\begin{align*}
\min_{0<{\text{\tt j}}<q_{n_k+1}}&\vvvert  Q{\text{\tt j}} \alpha\vvvert \geq \min \left( \vvvert q_{n_k}\alpha\vvvert ,\frac{1}{2Qq_{n_{k}+1}}\right)=\frac{1}{2Q q_{n_{k}+1}}\\ &\text{ implying }\text{\tt disc}(q_{n_k+1})\geq \frac{1}{2Q^2 q_{n_{k}+1}}.
\end{align*}

This proves (\Ladiesroom).\ \ \Checkedbox

\

\demo{Construction of measurable sets as in {\large\smiley}}

By   (\Ladiesroom) there exist a subsequence $(n_k)\uparrow\infty$ and $\theta>0$ such that $\text{\tt disc}(q_{n_k})>\frac{\theta}{q_{n_k}}$  and  such that $q_{n_k}$ is sufficiently large compared to $|F|^2$, 
where $F$ is the finite set of values taken by $\v_{q_{n_k}}$ as in the remark after the  Denjoy-Koksma inequality.

\ 

Fix $0\leq \epsilon\leq Q-1$. To obtain the periodicity $\Phi(\epsilon+1)-\Phi(\epsilon)$, we   build sequences of measurable sets $(A_k), (B_k)\subset \Bbb T$ such that
\bul $\v_{q_{n_k}}$ is constant on $A_k$ and $B_k$,
\bul $\v_{q_{n_k}}|_{B_k}-\v_{q_{n_k}}|_{A_k}=\Phi(\epsilon+1)-\Phi(\epsilon)$ and 
\bul $m_{\Bbb T}(A_k), m_{\Bbb T}(B_k)>c>0$. 

\

Fix $k$ and let $\bdy$ be the partition of $\Bbb T$  by the discontinuities $\{\frac{\epsilon}{Q}-h \alpha:\ 0\leq h\leq q_{n_k}-1\}$ of the step function $\v_{q_{n_k}}$.
\

For $0\le h<q_{n_k}$, let $I_h^-\in\bdy$ be the interval with right endpoint $\frac{\epsilon}{Q}-h \alpha$ and $I_h^+\in\bdy$ be the interval with left endpoint $\frac{\epsilon}{Q}-h \alpha$.
\

We can choose $0<h_1, h_2,\ldots, h_{\left\lfloor \frac{q_{n_k}}{|F|^2}\right\rfloor}<q_{n_k} $ such that $\v_{q_{n_k}}$ is constant on
$$A_k:=\bigcup_{u=1}^{\left\lfloor \frac{q_{n_k}}{|F|^2}\right\rfloor}I^-_{h_u}\text{ and }B_k=\bigcup_{u=1}^{\left\lfloor \frac{q_{n_k}}{|F|^2}\right\rfloor}I^+_{h_u}.$$

Evidently,
$$\v_{q_{n_k}}|_{B_k}-\v_{q_{n_k}}|_{A_k}=\Phi(\epsilon+1)-\Phi(\epsilon)$$
and by (\Ladiesroom)
$$m(A_k),m(B_k)\geq \text{\tt disc}(q_{n_k})\left\lfloor \frac{q_{n_k}}{|F|^2}\right\rfloor\geq \frac{\theta}{2|F|^2}.$$

      These sets are as in {\large\smiley} and the proof of theorem 1 in the countable case  is now complete.
      \demo{Proof of theorem 1 in the uncountable case}
      \
      
      Let 
      $$V:=\text{\tt Span}_{\Bbb Q}\Phi(\Bbb Z_Q)\subset\Bbb R^d,$$  let $K:=\dim\,V$ and let $\{e_k:\ 1\le k\le K\}$ be a basis for $V$ so that
      each
      $$\Phi(\e )=\sum_{k=1}^K\phi_k(\e )e_k\ \ \text{with}\ \phi_k(\e )\in\Bbb Z\ \ (1\le k\le K,\ \e \in\Bbb Z_Q).$$
      Consider the cocycle $\Psi:\Bbb T\to\Bbb Z^K$ defined by
      $$\Psi(x):=\phi(\lfl Qx\rfl)\ \text{where}\ \phi(\e ):=(\phi_1(\e ),\dots,\phi_K(\e ))\ \ (\e \in\Bbb Z_Q).$$ It follows that $\<\Psi(\Bbb T)\>=\Bbb Z^K$.
    We claim that
    $$\int_\Bbb T\Psi(x)dx=\tfrac1Q\sum_{\e \in\Bbb Z_Q}\phi(\e )=0.$$
    To see this,
      \begin{align*}0&=Q\int_\Bbb T\varphi(x)dx\\ &=\sum_{\e \in\Bbb Z_Q}\Phi(\e )\\ &=\sum_{k=1}^K\left(\sum_{\e \in\Bbb Z_Q}\phi_k(\e )\right)e_k.
      \end{align*}
By  linear independence of $\{e_k:\ 1\le k\le K\}$, for each $1\le k\le K$
$\sum_{\e \in\Bbb Z_Q}\phi_k(\e )=0$ showing that indeed
$\int_\Bbb T\Psi(x)dx=
0.$
\

Thus, by (\Wheelchair) as on page \pageref{wheelchair} in the countable case, and Schmidt's theorem,
$$\<\phi(\Bbb Z_Q)\>\subset\ \text{\tt Per}\,(\Psi)=E(\Psi).$$
It follows that $\Phi=L\circ\phi$ (and $\phi=L\circ\Psi$) where $L:\Bbb Z^K\to V\subset\Bbb R^d$ is given by
$$L(z_1,\dots,z_K):=\sum_{k=1}^Kz_ke_k.$$
By linearity of $L$,
$$L(E(\Psi))\subset E(L\circ\Psi)=E(\v)$$
and
$$\Phi(\Bbb Z_Q)=L(\phi(\Bbb Z_Q))\subset\ E(\v).\ \ \CheckedBox$$

\section*{\S2 The orbit sequence}
Theorems 2 and 3 both depend on the modeling of the {\it orbit sequence}
$$(\v_n(0):\ n\ge 1)$$ by an associated {\tt affine random walk}. To extract this affine random walk we first obtain a sequential 
substitution construction of the 
{\it jump sequence}
$$(\v(\{n\a\}):\  n\ge 1)\ \text{for}\ \ \a\in (0,1)\setminus\Bbb Q.$$

\

To this end, let $\b=\{Q\a\} \&\ P:=\lfl Q\a\rfl$ so that
$\a=\frac{P+\b}Q$.  
\

Define the map $\pi:[0,1)\to [0,1)\x\Bbb Z_Q$ by
$$\pi(x):=(\{Qx\} ,\lfl Qx\rfl),$$ 
the transformation $\tau:[0,1)\x\Bbb Z_Q\to[0,1)\x\Bbb Z_Q$ by $\tau:=\pi\circ r_\a\circ\pi^{-1}$ and $\text{\tt k}:[0,1)\times \Z_Q\to \Z_Q$ by
$$\text{\tt k}(x,k):=\upkappa\circ\pi^{-1}(x,k)=k;$$
then
\begin{align*}
 \tau(y,k)&=\pi(\{\tfrac{y+k}Q+\a\})\\ &=\pi(\{\tfrac{y+k+P+\b}Q\})\\ &=(r_{\b}(y),\lfl k+P+y+\b\rfl\ \mod Q)\\ &=(r_{\b}(y),k+P+1_{[1-\b,1)}(y)\ \mod Q).
\end{align*}

Thus
\begin{align*}\tag{\IroningI}\label{IroningI}
\upkappa(\{n\a\})&=\upkappa\circ r_\a^n(0)=\text{\tt k}\circ\pi\circ r_\a^n(0)\\ &=\text{\tt k}\circ\tau^n\pi(0)\\ &=nP+\sum_{k=0}^{n-1}1_{[1-\b,1)}(\{k\b\})\\ &=nP+\sum_{k=1}^{n}\psi_k\ \mod\ Q
\end{align*}
where $\psi_k:=1_{[1-\b,1)}(\{(k-1)\b\})$. The sequence $(\psi_k:\ k\ge 1)$ is generated as follows.

   The {\it modified continued fraction expansion} of $\b\in (0,1)$ is
   \begin{align*}\b & =\frac{1}{n_1-\frac{1}{n_2-_{\ddots-\frac{1}{n_k-_{\ddots}}}}}\\ & \,
\\ &=:\cfrac[r]{1|}{|n_1}-\cfrac[r]{1|}{|n_2}-\dots-\cfrac[r]{1|}{|n_k} \ \ -\cdots\\ &=
[n_1,n_2,\dots]\end{align*}
     where $n_k(\b):=n(s^{k-1}(\b))$ with $n(\b):=\lcl\frac1\b\rcl$ and
   $s(\b):=1-\{\frac1\b\}=n(\b)-\frac1\b$.

   See \cite{Keane1970}\ $\&$\ \cite{KN}.
   
    \subsection*{The quadratic case}\ \ If $\a\in\text{\tt QUAD}$, then so does $\b=\{Q\a\}$ and there exist
      $$(n_1,n_2,\dots,n_K)\in\Bbb N_2^K\ \&\ (m_1,\dots,m_L)\in\Bbb N_2^L\setminus \{2\mathbb{1}\}$$ such that
  $$\b=[n_1,n_2,\dots,n_K,\overline{m_1,\dots,m_L}].$$
  Here and throughout, 
 \bul $\mathbb{1}$ denotes a vector all of whose coordinates are $1$,\bul  $b_0\odot b_1$ denotes the concatenation of the finite sequences $b_0$ and $b_1$,
\bul  and 
$b_0^{\odot n}$  denotes  the concatenation of $n$ copies of $b_0$.
   \proclaim{Theorem 2.1 in \cite{AK}}
   For $\b=[n_1,n_2,\dots]$, let
   $b_0(0)=0,\ b_0(1)=1\ \&$
   $$b_{k+1}(0)=b_k(0)^{\odot (n_{k+1}-1)}\odot b_k(1)\ \ \&\ \ b_{k+1}(1)=b_k(0)^{\odot (n_{k+1}-2)}\odot b_k(1),$$
   then
   $$(\psi_1,\dots,\psi_{\ell_k(i)})=b_k(i)\ \ \ \ (k\ge 1)$$ if the last symbol in $b_k(1)$ is changed from ``1'' to ``0''.
   \endproclaim
  Here    $\ell_k(i)=|b_k(i)|\ \ (i=0,1)$ are the \emph{block lengths}.

   \subsection*{Block lengths}
   \
   
   Let $\u\ell_k:=\begin{pmatrix}
\ell_k(0)\\ \ell_k(1)
\end{pmatrix}$, then
$$\u\ell_0=\begin{pmatrix}
1\\ 1
\end{pmatrix}\ \&\ \u\ell_{k+1}=\begin{pmatrix} n_{k+1}-1 & 1\\
n_{k+1}-2 & 1
\end{pmatrix}\u\ell_k.$$
   \subsection*{Parities, Jumps $\&$  Orbits}
   Next, we compute the {\tt jump} blocks. 
   \
   
   We call $\upkappa(x)=\lfl Qx\rfl\in\Bbb Z_Q$ the {\it parity} of $x$ and we begin
   by calculating the {\it parity blocks} with a generalization of \cite[theorem 2.2]{AK}.
\

   For $\b=\{Q\a\}=[n_1,n_2,\dots],\ \e\in\Bbb Z_Q,\ i=0,1$,  define
      \begin{align*}
   B_k(i,\e)&:=(\e+\upkappa(\{(n-1)\a\}): 1\le n\le \ell_k(i)),\end{align*}
   then by (\IroningI) as on page \pageref{IroningI} and \cite[theorem 2.1]{AK} respectively,
   
   \begin{align*}\label{Bat}\tag{\Bat} B_k(i,\e)&= (\e+(n-1)P+\sum_{\nu=1}^{n-1}1_{[1-\b,1)}(\{(\nu-1)\b\}):\ 1\le n\le \ell_k(i))\\ &=
   (\e+(n-1)P+\sum_{\nu=1}^{n-1}b_k(i)_\nu:\ 1\le n\le \ell_k(i))\end{align*}
   where the addition is $\ \mod Q$ and $\sum_{\nu\in\emptyset}:=0$. Note that $B_0(i, \epsilon)=(\epsilon)$.
   \proclaim{Theorem 5.1 \ \ (Parity recursions)}
    \begin{align*}\tag{\Coffeecup} \label{Coffeecup}
        B_{k+1}&(i,\e)=\\ &\bigodot_{j=1}^{n_{k+1}-1-i}B_k(0,\e+(j-1)\e_k)\odot B_k(1,\e+(n_{k+1}-1-i)\e_k).\end{align*}
  with $\e_k:=\sum_{j=1}^{\ell_k(0)}(b_k(0))_j+\ell_k(0)P\ \mod Q$ and $\bigodot_{j\in\emptyset}H_j\odot B:=B$ for finite sequences $(H_j)$ and $B$. Here (as before) the addition is $\mod\ Q$.\endproclaim
  
It follows that 
$$B_1(i,\e)=(\e,\e+P,\dots,\e+(n_1-1-i)P) \mod Q.$$
  
\demo{Proof}  Fix $i=0,1,\ \e\in\Bbb Z_Q, k\ge 1$ and $1\le  n\le\ell_{k+1}(i)$, then $n=q\ell_k(0)+r$ where $0\le q\le n_{k+1}-i-1$ and $1\le r\le\ell_k(j_q)$ with
$j_{n_{k+1}-i-1}=1\ \&\ j_q=0$ for $q<n_{k+1}-i-1$.

Using \cite[theorem 2.1]{AK}  and (\Bat) as on page \pageref{Bat}, we have $\mod Q$,
\begin{align*}B_{k+1}(i,\e)_n &=\e+(n-1)P+\sum_{{\nu}=1}^{n-1}b_{k+1}(i)_{\nu}\\ &=
\e+(q\ell_k(0)+r -1)P+\sum_{{\nu}=1}^{q\ell_k(0)+r -1}b_{k+1}(i)_{\nu}\\ &= \e+q\e_k+(r -1)P+\sum_{{\nu}=1}^{r -1}b_{k}(j_q)_{\nu}\\ &=
\(\bigodot_{j=1}^{n_{k+1}-1-i}B_k(0,\e+(j-1)\e_k)\odot B_k(1,\e+(n_{k+1}-1-i)\e_k)\)_n.\ \ \CheckedBox
\end{align*}

\subsection*{Parity states and transition algorithm}
\

The $k^{\text{\tiny th}}$ {\it parity states} are $\e_k(i)\ \ (i=0,1)$ where $$\e_k(i):=\sum_{j=1}^{\ell_k(i)}(b_k(i))_j+\ell_k(i)P\ \mod Q.$$
In (\Coffeecup) as on page \pageref{Coffeecup}, $\e_k=\e_k(0)$.

\

The parity states are given by $\e_0(i)=P+i\mod Q$ and 
\begin{align*}\tag{\phone}\label{phone}  \e_{k+1}(i)=(n_{k+1}-i-1)\e_k(0)+\e_k(1)\mod Q\  \ (i=0,1,\ k\ge 1). 
\end{align*}

\proclaim{Parity proposition}\ \ For every $k\ge 1$, $\<\{\e_k,\e_{k+1}\}\>=\Bbb Z_Q$.\endproclaim\demo{Proof}
Define $\zeta_k=(\zeta_k(0),\zeta_k(1))$ by
$$\zeta_0(i)=P+i\ \&\ \zeta_{k+1}(i)=(n_{k+1}-i)\zeta_k(0)+\zeta_k(1).$$
It follows that
\bul\ $\e_k(i)=\zeta_k(i)\ \mod Q$;
\bul $\text{\tt gcd}\,\{\zeta_k(0),\zeta_k(1)\}=1\ \forall\ k\ge 1$;
\bul $\text{\tt gcd}\,\{\zeta_k(0),\zeta_{k+1}(0)\}=1\ \forall\ k\ge 1$;
\bul $\<\{\zeta_k(0),\zeta_{k+1}(0)\}\>=\Bbb Z\ \forall\ k\ge 1$;
\bul $\<\{\e_k,\e_{k+1}\}\>=\Bbb Z_Q$.\ \Checkedbox
  
   \subsection*{Jump blocks}
   \
   
Next, for $k\ge 1,\ \e\in\Bbb Z_Q\ \&\ i=0,1$, define the {\it auxiliary jump blocks}
   $$J_k(i,\e):=\Phi(B_k(i,\e))$$
   where
   $$\Phi((a_1,\dots,a_n)):=(\Phi(a_1),\dots,\Phi(a_n)).$$
      \
   
   It follows from (\Coffeecup) that for $i=0,1$:
   \begin{align*}\tag{*}
         J_{k+1}(i,\e)=\bigodot_{j=1}^{n_{k+1}-1-i}J_k(0,\e+(j-1)\e_k)\odot J_k(1,\e+(n_{k+1}-1-i)\e_k);
 \end{align*}
  where addition is $\mod Q$ and that the {\it jump block} 
      \begin{align*} \tag{\sun} (\v(\{j\a\}))_{j=0}^{\ell_k(0)-1}=J_k(0,0).     
      \end{align*}

  \subsection*{  Orbit blocks}\ \ 
  \

  Define the {\it  auxiliary orbit blocks}
 $$\Si_k(i,\e):=(\sum_{{\nu}=1}^jJ_k(i,\e)_{\nu}:\ 1\le j\le \ell_k(i)).$$
  In particular by  (\sun)
   $$\Si_k(0,0)=(\v_1(0),\v_2(0),\dots,\v_{\ell_k(0)}(0)).$$

Our goal here is to obtain the transition between auxiliary orbit blocks.
\subsection*{  Orbit block  transitions} \
 
 \
 
        The {\it simple displacement} over the auxiliary jump block $J_k(i,\e)$ is 
        $$\s_k(i,\e):=\Si_k(i,\e)_{\ell_k(i)}=\sum_{j=1}^{\ell_k(i)}J_k(i,\e)_j.$$
        
      The  {\it cumulative displacements} over the concatenation jump blocks $\bigodot_{j=1}^K J_k(0,\e+(j-1)\e_k \mod Q)\ \ (K\ge 0)$ are
  $$s_k(K,\e):=\sum_{j=1}^K\s_k(0,\e+(j-1)\e_k\mod Q).$$      
By (*), for $k\ge 1,\ \e\in\Bbb Z_Q,\ i=0,1$, 
   \begin{alignat*}{2}&\ \ \ \ \Si_{k+1}(i,\e)=\bigodot_{j=1}^{n_{k+1}-1-i}(\Si_k&&(0,\e+(j-1)\e_k\mod Q)+s_k(j-1,\e)\mathbb{1})\odot\\ &  &&
   \odot(\Si_k(1,\e+(n_{k+1}-i-1)\e_k\mod Q)+s_k(n_{k+1}-i-1,\e)\mathbb{1}).
   \end{alignat*}

   \subsection*{ Generating functions of orbit blocks}
    \
    
    For $k\geq 1$ define the functions
  $x_k(i,\e):\Om_k(i)=[1,\ell_k(i)]\to\Bbb R^d$ by
  \begin{align*}\tag{\Biohazard}\label{Biohazard}
     x_k(i,\e)(\om):= \Si_{k}(i,\e)_\om\ \ \ (\om\in\Om_k(i))
  \end{align*}

  and their {\it generating functions}
  $$U_k(i,\e,\th):=\sum_{\om\in\Om_k(i)}e^{2\pi\frak i \<\th, x_k(i,\e)(\om)\>}\ \ \ \ (\th\in\Bbb T^d).$$
  Here and throughout, $\frak i:=\sqrt{-1}$.
\subsection*{Transition matrices}
  Noting that
  $$\Om_{k+1}(i)=\odot_{j=1}^{n_{k+1}-i-1}(\Om_k(0)+(j-1)\ell_k(0))\odot (\Om_k(1)+(n_{k+1}-i-1)\ell_k(0)),$$
  we have
  {\footnotesize\begin{align*}&U_{k+1}(i,\e,\th)=\sum_{\om\in\Om_{k+1}(i)}e^{2\pi\frak i \<\th,x_{k+1}(i,\e)(\om)\>}\\ &=
  (\sum_{j=1}^{n_{k+1}-i-1}\sum_{\om\in\Om_k(0)+(j-1)\ell_k(0)}+\sum_{\om\in \Om_k(1)+(n_{k+1}-i-1)\ell_k(0)})e^{2\pi\frak i \<\th,x_{k+1}(i,\e)(\om)\>}\\ &=
  \sum_{j=1}^{n_{k+1}-i-1}\sum_{\om\in\Om_k(0)}e^{2\pi\frak i \<\th, (x_k(0,\e+(j-1)\e_k)(\om)+s_k(j-1,\e))\>}+
  \sum_{\om\in \Om_k(1)}e^{2\pi\frak i \<\th, (x_k(1,\e+(n_{k+1}-i-1)\e_k)(\om)+s_k(n_{k+1}-i-1,\e))\>}\\ &=
  \sum_{j=1}^{n_{k+1}-i-1}e^{2\pi\frak i \<\th, s_k(j-1,\e)\>}U_k(0,\e+(j-1)\e_k,\th)+e^{2\pi\frak i \<\th, s_k(n_{k+1}-i-1,\e)\>}U_k(1,\e+(n_{k+1}-i-1)\e_k,\th)\\ &=
  \sum_{\D\in\Bbb Z_Q}\sum_{j\in\frak m(\e_k,\D)\cap [1,n_{k+1}-i-1]}e^{2\pi\frak i \<\th, s_k(j-1,\e)\>}U_k(0,\e+\D,\th)+e^{2\pi\frak i \<\th, s_k(n_{k+1}-i-1,\e)\>}U_k(1,\e+(n_{k+1}-i-1)\e_k,\th)
  \end{align*}}
   where  for $\e,\ \D\in\Bbb Z_Q$, 
   $$\frak m(\e,\D):=\{j\in\Bbb N:\ (j-1)\e=\D\ \mod\  Q\}$$ (with  $\sum_{\om\in\emptyset}:=0$ as before).

Equivalently,
  $$U_{k+1}(\th)=A^{(k+1)}(\th)U_k(\th)$$
  where $S:=\{0,1\}\x\Bbb Z_Q$ and $U_k: \mathbb T^d \to \mathbb C^S$ is given by
  \begin{align*}
     U_k(\theta)_{(i, \epsilon)}&:=U_k(i,\e, \theta)  \text{ for } (i,\e)\in S\ \ \&\\ & A^{(k+1)}:\Bbb T^d\to M_{S\x S}(\Bbb C):=\{a:S\x S\to \Bbb C\}
  \end{align*}
              
is given by:

 \begin{align*}
& A^{(k+1)}_{(i,\e),(0,\e+\D)}(\th)=\sum_{j\in \frak m(\e_k,\D)\cap [1,n_{k+1}-i-1]}e^{2\pi\frak i \<\th, s_k(j-1,\e)\>}\ \ \ \ \text{if}\ \ (n_{k+1},i)\ne (2,1),
\\ & A^{(k+1)}_{(i,\e),(0,\e+\D)}(\th)=0  \ \ \ \ \text{if}\ \ (n_{k+1},i)= (2,1),\\ &                     
   A^{(k+1)}_{(i,\e),(1,\e+\D)}(\th)=e^{2\pi\frak i \<\th, s_k(n_{k+1}-i-1,\e)\>}1_{\{\D\}}((n_{k+1}-i-1)\e_k).
   \end{align*}

  \
  
  It follows  that
  \begin{align*}A^{(k+1)}_{(i,\e),(0,\e+\D)}(0)=N_{k+1}(i,\D)\ \&\ \ A^{(k+1)}_{(i,\e),(1,\e+\D)}(0)=1_{\{\D\}}((n_{k+1}-i-1)\e_k)  
  \end{align*}

where   $N_{k+1}(i,\D):=\#\,\frak m(\e_k,\D)\cap [1,n_{k+1}-i-1]$.

    \section*{\S3  The random affine model}

 \subsection*{Probabilities}
Here, we consider the probabilities
$$P_k^{(i)}:=\frac{\#}{\ell_k(i)}\in\mathcal P(\Om_k(i))$$ and each $x_k(i,\e)\in\Bbb R^d$ as a random variable with sample space
$(\Om_k(i),P_k^{(i)})$, and understand the transitions of the resulting stochastic processes 
$$(x_k(i,\e):\  k\ge 1)\ \ ((i,\e)\in S)$$ 
 in the {\tt RAT} lemma.

Let
$$\Xi_k(i,\e,\th):=E(e^{2\pi\frak i \<\th, x_k(i,\e)\>})=\frac1{\ell_k(i)}U_k(i,\e,\th),$$ then
 \begin{align*}
  \tag{\dsjuridical}\label{jur}\Xi_{k+1}(\th)=\Pi^{(k+1)}(\th)\Xi_k(\th)
 \end{align*}
where $\Xi_k:=(\Xi_k(i,\e):\ \ (i,\e)\in S)\ \ \&\ \Pi^{(k+1)}(\th):S\x S\to\Bbb C$ is given by
 \begin{align*}
\Pi^{(k+1)}_{(i,\e),(j,\D)}(\th)=\frac{\ell_k(j)}{\ell_{k+1}(i)}A^{(k+1)}_{(i,\e),(j,\D)}(\th). \end{align*}

 \
 
 \subsection*{Random variables}
 \
 
We denote by $\text{\tt RV}\,(Z)$, for $Z$  a measurable space the collection of $Z$- valued random variables.
 
Consider any sequence of independent,  random vectors 
 \begin{align*}\tag{\Dontwash}\label{Dontwash}
  (\mathcal L^{(k+1)}_s,\ W_{s,t}^{(k+1)}:\ \ s,t\in S)\in \text{\tt RV}\,(S^S\x (\Bbb R^d)^{S\x S})\ \ (k\ge 0)
 \end{align*}
whose marginals satisfy
    \begin{align*}&P(\mathcal L^{(k+1)}_{(i,\e)}=(0,\e+\D))=\frac{\ell_k(0)N_{k+1}(i,\D)}{\ell_{k+1}(i)},\\ &\ 
    P(\mathcal L^{(k+1)}_{(i,\e)}=(1,\e+\D))=\frac{\ell_k(1)}{\ell_{k+1}(i)}1_{\{\D\}}((n_{k+1}-i-1)\e_k);
    \\ &P([W^{(k+1)}_{(i,\e),(0,\e+\D)}=s_k(J-1,\e)]|[\mathcal L^{(k+1)}_{(i,\e)}=(0,\e+\D)])=\\  
&=\frac{\#\{j\in\frak m(\e_k,\D)\cap [1,n_{k+1}-i-1]:\ s_k(j-1,\e)=s_k(J-1,\e)\}}{N_{k+1}(i,\D)};\\ &\text{\tt for }  \ J\in\frak m(\e_k,\D)\cap [1,n_{k+1}-i-1]\ \&\\ &
P([W^{(k+1)}_{(i,\e),(1,\e+\D)}=s_k(n_{k+1}-i-1,\e)]|[\mathcal L^{(k+1)}_{(i,\e)}=(1,\e+\D)])=1.
      \end{align*}
     
      Note that  when $n_{k+1}=2$, then  $\mathcal L^{(k+1)}_{(1,\e)}=(1,\e+(n_{k+1}-i-1)\e_k)$ a.s. $\forall\ \e\in\Bbb Z_Q$
         and that   
   $W_{s,t}^{(k+1)}$ is defined when and only when $P(\mathcal L^{(k+1)}_s=t)>0$.
   
   \

 \subsection*{Random affine transformations }
 \
 
 Given a random vector
  $$(\mathcal L,W)\in\text{\tt RV}(S^S\x (\mathbb R^d)^{S\x S}),$$ the associated
  {\it random affine transformation} ({\tt RAT}) $F\in\text{\tt RV}\,(M_{S\x S}(\Bbb Z)\x (\mathbb R^d)^S)$ defined by
  \begin{align*}\tag{\Handwash}
     F(x)_s:=x_{\mathcal L_s}+W_{s,\mathcal L_s}=:(a(F)x)_s+b(F)_s\text{ for }x\in (\R^d)^S.
  \end{align*}

This {\tt RAT} is of {\it flip-type} in the sense of \cite{ABN2016}.
\

Throughout this paper we'll often denote a flip-type {\tt RAT} $$F=(a(F),b(F))\in\text{\tt RV}\,(M_{S\x S}(\{0,1\})\x(\Bbb R^d)^S)$$ by
 $$F=(\mathcal L,W)=(\mathcal L(F),W(F))\in\text{\tt RV}\,(S^S\x (\Bbb R^d)^{S\x S}).$$ Here 
 $$a_{s,t}=\d_{t,\mathcal L_s}\ \&\ b_s=W_{s,\mathcal L_s}.$$

 Given a  sequence $(\mathcal L^{(k+1)}_s,\ W_{s,t}^{(k+1)}:\ \ s,t\in S)\ \ (k\ge 0)$ of independent random vectors 
 as before, consider  the associated {\tt RAT} sequence $$(F_k:\ k\ge 1)\in \text{\tt RV}\,(M_{S\x S}(\{0,1\})\x(\Bbb R^d)^S)^\Bbb N$$ of
 independent {\tt RATs} defined by (\Handwash).

  \subsection*{{\tt RAT} characteristic function}
  \
  
The  {\it characteristic function} of the  {\tt RAT} $F=(\mathcal L,W)\in\text{\tt RV}\,(S^S\x (\Bbb R^d)^{S\x S})$ ( {\tt RAT-CF}) 
  is $\Pi_F:\Bbb R^d\to M_{S\x S}(\Bbb C)$ defined by
\begin{align*}\tag{\dsliterary}\label{lit}\Pi_F(\th)_{s,t}=P(\mathcal L_s=t)E(e^{2\pi\frak i\<\th,W_{s,t}\>})\ \ \ (s,t\in S). 
\end{align*}
Note that $\Pi^{(k+1)}$ in (\dsjuridical)  on page\ \ \pageref{jur}\  \ is the {\tt RAT-CF} of the {\tt RAT} 
$(\mathcal L^{(k+1)},W^{(k+1)})$ where $\mathcal L^{(k+1)}\ \&\ W^{(k+1)}$ are as in (\Dontwash)  on page \ \pageref{Dontwash}.
  \proclaim{{\tt RAT} lemma }
 \
 
 \ \ For each $k\ge 1,\ s\in S$:
 $$\text{\tt dist}\,F_1^k(\u 0)_s\ =\ \text{\tt dist}\,x_k(s)$$ where $x_k(s)$ is as in (\Biohazard) as on \pageref{Biohazard}.\endproclaim
Here and throughout for $K\le L\ \&$ {\tt RAT}s $(F_j:\ K\le j\le L)$
 $$F_K^L:=F_L\circ F_{L-1}\circ \dots\circ F_{K+1}\circ F_K.$$\demo{Proof}\ \ For $k\ge 1$, define
 $$X^{(k)}:=F_k\circ F_{k-1}\circ\dots\circ F_1(0)$$ and
 $$\widehat{\Xi}_k:=(E(e^{2\pi\frak i \<\th, X^{(k)}(i,\e)\>}):\ \ (i,\e)\in S).$$
 By construction,
 \begin{align*}
\widehat{\Xi}_{k+1}(\th)=\Pi^{(k+1)}(\th)\widehat{\Xi}_k(\th).
 \end{align*}
 By (\dsjuridical) as on page \pageref{jur},\begin{align*}
\Xi_{k+1}(\th)=\Pi^{(k+1)}(\th)\Xi_k(\th).
 \end{align*}
 
 The result follows by induction since $\Xi_0=\widehat{\Xi}_0\equiv \mathbb{1}$.\ \ \Checkedbox
 \
 
   \subsection*{Associated affine random walks}
   \

We associate to  a sequence $$(F_k:\ k\ge 1)\in \text{\tt RV}\,(M_{S\x S}(\{0,1\})\x(\Bbb R^d)^S)^\Bbb N$$ of
 independent {\tt RATs} an   {\it affine random walk} ({\tt ARW}).
 \

 This is the  $(\Bbb R^d)^S$-valued stochastic process 
 $$(X^{(k)}=(X^{(k)}_s:\ s\in S):\ k\ge 1)$$
defined by $$X^{(k)}:=F_1^k(0).$$
\

 \subsection*{Elementary presentation\label{elpres}}
 \

   We now split the random  vectors $(\mathcal L^{(k+1)}_s,\ W_{s,t}^{(k+1)}:\ \ s,t\in S)\ \ (k\ge 0)$    
   into more elementary components.

   Write
$$\mathcal L^{(k+1)}_{(i,\e)}=:(\frak r^{(k+1)}_{(i,\e)},\frak s^{(k+1)}_{(i,\e)}),$$
then
 $\frak r^{(k+1)}_{(i,\e)}=\frak r^{(k+1)}_i$ is a $\{0,1\}$-valued random variable where
   $$P(\frak r^{(k+1)}_i=0)=\frac{\ell_k(0)(n_{k+1}-i-1)}{\ell_{k+1}(i)}\ \&\ P(\frak r^{(k+1)}_i=1)=\frac{\ell_k(1)}{\ell_{k+1}(i)}$$
   and
   $$\frak s^{(k+1)}_{(i,\e)}=\e+\frak e^{(k+1)}_i\mod Q$$
where $\frak s^{(k+1)}_{(i,\e)}$ and $\frak e^{(k+1)}_i$are $\Z_Q$-valued random variables; the latter is given by
 $$P(\frak e^{(k+1)}_i=\D\|\frak r^{(k+1)}_i=0)=\frac{N_{k+1}(i,\D)}{n_{k+1}-i-1}\ \ \ (\D\in\Bbb Z_Q)$$
 and
  $$P(\frak e^{(k+1)}_i=(n_{k+1}-i-1)\e_k\|\frak r^{(k+1)}_i=1)=1.$$
  Next define random variables $\frak u^{(k+1)}(i)\ \ (k\ge 1,\ i=0,1)$ by
  $$\frak u^{(k+1)}(i)\ \ \begin{cases} & \text{uniform on $\frak m(\e_k,\frak e^{(k+1)}_i)\cap [1,n_{k+1}-i-1]$ if $\frak r_i^{(k+1)}=0$}\\ & \equiv\ 
                           n_{k+1}-i\ \ \ \text{if}\ \ \ \frak r_i^{(k+1)}=1.
                          \end{cases}$$

  Now we define  random variables $W^{(k)}_s\ \ (k\ge 1,\ s\in S)$ by
  $$W_{(i,\e)}^{(k+1)}:=s_k(\frak u^{(k+1)}(i)-1,\e).$$
  It is not hard to see that
  $$P(W^{(k+1)}_{(i,\e),(j,\e+\D)}=J\|\mathcal L^{(k+1)}_{(i,\e)}=(j,\e+\D))=P(W_{(i,\e)}^{(k+1)}=J\|\frak r^{(k+1)}_{i}=j,\ \frak e^{(k+1)}_i=\D).$$
  In the sequel, we'll have recourse to the {\it elementary random vector sequence} $(\frak x^{(k)}:\ k\ge 1)\in RV((\{0,1\}\x\Bbb Z_Q\x\Bbb N_0)^{\{0,1\}})^\Bbb N$ where
  $$\frak x^{(k)}:=(\frak r^{(k)}_i,\ \frak e^{(k)}_i,\ \frak u^{(k)}(i):\ i=0,1).$$
  The {\tt RAT} $F_k$ is constructed from $\frak x^{(k)}\ \&$ the (deterministic) cumulative displacements $s_{k-1}$.

 \

\section*{\S4 The {\tt RAT} sequence in the quadratic case}
   \

  We assume that $\alpha\in \text{\tt QUAD}$; thus $\b=\{Q \alpha\}\in\text{\tt QUAD}$. These hold if and only if there exist
      $$(n_1,n_2,\dots,n_K)\in\Bbb N_2^K\ \&\ (m_1,\dots,m_L)\in\Bbb N_2^L\setminus 2\mathbb{1}$$ such that
    \begin{align*}
        \tag*{\Football}\label{Football}\beta:=[n_1,n_2,\dots]=:[n_1,\dots,n_K,\overline{m_1,\dots,m_L}].
    \end{align*}
\

We next establish that the centered {\tt RAT} sequence (as in \cite{ABN2016}) corresponding to a quadratic irrational and a rational step function is 
``{\tt asymptotically eventually periodic}''.
\

The proofs of theorems 2 $\&$ 3 rely on this fact.

\

This asymptotic eventual periodicity is obtained via centering.  We'll see that elementary random vector sequence is always  asymptotically eventually periodic,
however, the cumulative displacements may have linear growth. The centering is needed to offset this possibility.

\

In this section, we'll often ``{\tt possibly extend the period in \Football}'' to demonstrate 
eventual periodicity of related sequences.

\

This means that for some $M\in\Bbb N$, we'll modify \Football\ to
$$[n_1,n_2,\dots]=[n_1,\dots,n_K,\overline{\undersetbrace{\text{\tiny $M$-times}}{m_1,\dots,m_L,\dots,m_1,\dots,m_L}}].$$
Recall ((\phone) on page \pageref{phone})  that the parity state transitions are given by
  $$\e_{k+1}(i)=(n_{k+1}-i-1)\e_k(0)+\e_k(1)\ \ \mod\ Q.$$
 In the quadratic case, these transitions form an  eventually periodic sequence, whence 
 $$((\e_k(0),\e_k(1)):\ \ k\ge 1)$$ is also eventually periodic.
 \
 
  These parity transitions only depend on $\a\in\Bbb T\setminus\Bbb Q\ \&\ Q\ge 2$.
  \

 If $\a\in\text{\tt QUAD}$, then by possibly extending the period in \Football, we may assume that $(\e_{k+L}(0), \e_{k+L}(1))=(\e_k(0),\e_k(1))\ \forall\ k>K.$
   \subsection*{Simple displacement  transitions}
   \
   
  Consider the {\it simple displacement  vectors} $$\s_k:=(\s_k(i,\e):\ (i,\e)\in S)\in(\Bbb R^d)^S.$$
   By theorem 5.1,
   for $i=0,1$ and $(n_{k+1},i)\ne (2,1)$:
   \begin{align*}
         \s_{k+1}(i,\e)=\sum_{j=1}^{n_{k+1}-1-i}\s_k(0,\e+(j-1)\e_k)+\s_k(1,\e+(n_{k+1}-1-i)\e_k);
 \end{align*}
   where 
   $\e_k=\e_k(0)$ as before.
  
  Thus there exist matrices $M^{(k+1)}:S\times S\to \Z$ such that 
  $$\s_{k+1}=M^{(k+1)}\s_{k}$$
  for each $k\geq 1$.
  \

  The simple displacement  transformations also only depend on $\a\in\Bbb T\setminus\Bbb Q\ \&\ Q\ge 2$. 
  \
  
 Seeing $\s_k=(\s_k(s):\ s\in S)\in (\Bbb R^d)^S$ as 
 $$\s_k=((\s^{(j)}_k(s):\ s\in S):\ 1\le j\le d)\in (\Bbb R^S)^d,$$ we note that
 each $\s^{(j)}_k\in\Bbb R^S$ is a linear image of $\Phi^{(j)}\in \Bbb R^Q$ (the $j^{\text{\tiny th}}$ coordinate of $\Phi$) and
 $\s^{(j)}_{k+1}=M^{(k+1)}\s^{(j)}_{k}$ for each $1\leq j\leq d$.
  \

   \
   
\proclaim{Displacement lemma} 
\ Suppose that $\a\in\text{\tt QUAD}$, then there exist $K,\ L\in\Bbb N$ and $\frak c,\ \frak d\in(\Bbb R^d)^S$ so that
     \begin{align*}
     \s_{K+Ln}=\frak c+n\frak d.  \end{align*}
\endproclaim
For $\a\in\text{\tt QUAD}$, the simple displacement  transitions are eventually periodic and 
 the proof of the displacement lemma rests on the Denjoy-Koksma inequality and 
 a spectral analysis of the simple displacement  transformations on $\Bbb C^S$ over a period (as in the ``eigenvalue lemma'' below).

 {\subsection*{Subspace decomposition $\&$ eigenvalues}
 \
 
 For $\a\in\text{\tt QUAD}$, the parity sequence $(\e_k:\ k\ge 1)$ is eventually periodic, whence
 the above sequence of matrices $(M^{(k)}:\ k\ge 1)$ giving the displacement transitions is  also   eventually periodic. 
 \
 
 Suppose that
  \begin{align*}
    &[n_1,n_2,\dots]=[n_1,\dots,n_K,\overline{m_1,\dots,m_L}];\\ & (\e_k:\ k\ge 1)=(\e_1,\dots,\e_K,\overline{\eta_1,\dots,\eta_L});\\ &
   (M^{(k)}:\ k\ge 1)=(M^{(1)},\dots,M^{(K)},\overline{E_1,\dots,E_L}).
  \end{align*}

 Thus
 $$\s_{K+Ln}=B^n\s_K\ \text{where}\ B=E_L\cdots E_1.$$ 
 
 Next, write $\Bbb C^S=(\Bbb C^Q)^{\{0,1\}}$ and $z\in\Bbb C^S$ as $z=(z^{(0)},z^{(1)})\in (\Bbb C^Q)^{\{0,1\}}$.
 \
 
 The parity state transitions can now be rewritten as
 $$\s_{k+1}=M^{(k+1)}\begin{pmatrix}
\s{(0)}_k\\ \s{(1)}_k
\end{pmatrix}$$
where

$$\s{(i)}_k(\e)=\s_k(i,\e)\text{ for }\epsilon\in \Z_Q$$
and
\begin{align*}\tag{\Pointinghand}\label{Pointinghand}
M^{(k+1)}&=\frak P^{(k+1)}(\rho_{\epsilon_k})\\ &=\begin{pmatrix}
\frak P^{(k+1)}_{(0,0)}(\rho_{\epsilon_k})&\frak P^{(k+1)}_{(0,1)}(\rho_{\epsilon_k})\\
\frak P^{(k+1)}_{(1,0)}(\rho_{\epsilon_k})&\frak P^{(k+1)}_{(1,1)}(\rho_{\epsilon_k})
\end{pmatrix}\\ &
=\begin{pmatrix}
p_{n_{k+1}}(\rho_{\epsilon_k})&q_{n_{k+1}}(\rho_{\epsilon_k})\\
p_{n_{k+1}-1}(\rho_{\epsilon_k})&q_{n_{k+1}-1}(\rho_{\epsilon_k})
\end{pmatrix}
\end{align*}
with 
    $\rho_\epsilon\in M_{\Z_Q\times \Z_Q}(\Bbb C)$ defined by
    $$\rho_\e\, z(\d):=z(\d+\e)$$
and 
$$\frak P^{(k+1)}_{(0,0)},\frak P^{(k+1)}_{(0,1)},\frak P^{(k+1)}_{(1,0)},\frak P^{(k+1)}_{(1,1)},  p_\nu, q_\nu$$ 
are polynomials given by 
\begin{eqnarray*}
 p_\nu(x)&:=&\sum_{j=1}^{\nu-1}x^{j-1},\\
q_\nu(x)&:=&x^{\nu-1}\text{ and }\\
\frak P^{(k+1)}(x):=
\begin{pmatrix}
\frak P^{(k+1)}_{(0,0)}(x)&\frak P^{(k+1)}_{(0,1)}(x)\\
\frak P^{(k+1)}_{(1,0)}(x)&\frak P^{(k+1)}_{(1,1)}(x)
\end{pmatrix}
&:=&\begin{pmatrix}
p_{n_{k+1}}(x)&q_{n_{k+1}}(x)\\
p_{n_{k+1}-1}(x)&q_{n_{k+1}-1}(x)
\end{pmatrix}.
\end{eqnarray*}

Set $\g_r=e^{\frac{2 \pi i r}{Q}}$
and let $\vec e_{r}\in \mathbb C^Q$ be given by
$$(\vec e_{r})_{s}:=\g_{rs}$$
for $0 \leq r\leq Q-1$ and $1\leq s\leq Q$. 
\

Since  $\vec e_{s}\perp\vec e_{t}\ \forall\ s,\ t\in\Bbb Z_Q,\ s\ne t$, we have that
 $(\vec e_r:\ 0\le r\le Q-1)$ form an orthogonal  basis for $\mathbb C^Q$ and
\

\begin{equation*}
\label{equation: basis for null space} \text{\tt Span}\,\{\vec e_r~:~1\leq r \leq Q-1\}=\mathbb{1}^\perp=:\{\vec v=(v_h)\in \mathbb C^{Q}~:~ \sum_{h=0}^{Q-1}v_h=0\}.
\end{equation*}
Moreover,
\begin{equation*}\tag{\dstechnical}
 \rho_\e\vec e_r=\g_{r\e}\vec e_r.
\end{equation*}

Next, define the bracket $[\cdot,\cdot]:\Bbb C^{\{0,1\}}\x\Bbb C^Q\to \Bbb C^S=(\Bbb C^Q)^{\{0,1\}}$ by
$$[\vec c,\vec z]:=\begin{pmatrix}
c_0\vec z\\c_1\vec z
\end{pmatrix}$$ where $\vec c=(c_0,c_1)$.
\

It follows from (\dstechnical) that
$$M^{(k+1)}[\vec c,\vec e_r]=\frak P^{(k+1)}(\rho_{\epsilon_k})[\vec c,\vec e_r]=[\frak P^{(k+1)}(\g_{\epsilon_k r})\vec c,\vec e_r].$$

To summarize, letting  for $0\le r\le Q-1$,
$$V_{r}:=\{[\vec c,\vec e_r]:\ \vec c \in \mathbb C^{\{0,1\}}\},$$ then 
\sms  $\bigoplus_{r=0}^{Q-1}V_r=(\Bbb C^{Q})^{\{0,1\}}$ and   $BV_r=V_r\ \ (0\le r\le Q-1)$.
\proclaim{Eigenvalue lemma}

\sms For $1\le r\le Q-1$, all the eigenvalues of $B|_{V_r}$ are roots of unity.\endproclaim
\demo{Proof}\ \ 
\

We have that $B|_{V_0}$  is a product of integer matrices of the form 
$\begin{pmatrix} N& 1\\ N-1&1\end{pmatrix}$ with $N\in\Bbb N$; we have $N\ge 2$ for at least  one of these matrices. Therefore $B|_{V_0}$ is a positive matrix with integer coefficients and unit
determinant. 
It follows that the characteristic polynomial of $B|_{V_0}$ is an integer polynomial of form $z^2-Jz+1$ for some $J\in\Bbb N$ (and that $B|_{V_0}$  is hyperbolic). 
\

For each $1\le r\le Q-1$, 
$$|\det B|_{V_r}|=|\det\frak P^{(k+1)}(\g_{\epsilon_k r})|=1.$$
\

We claim first  that no $B|_{V_r}\ \ (1\le r\le Q-1)$ is hyperbolic. 
If this were not the case for $1\le r\le Q-1$, there would be $\l>1$ and a rational cocycle $ \Phi(\neq 0)\perp\mathbb{1}$ with $\<\Phi,\vec e_r\>\ne 0$ giving rise
to either 
\bul $\|\s_{K+L{n}}\|\gg \l^{n}$ which is impossible by the Denjoy-Koksma estimate;
\

or
\bul $\|\s_{K+L{n}}\|\ll  \frac1{\l^{n}}$ which is impossible by theorem 1.

\

To continue, since $B$ is an integer matrix, $\det (B-z\text{\tt Id})$ is a polynomial with integer coefficients. 
\

 It follows that
$$\det (B-z\text{\tt Id})|_{V_0^\perp}=\frac{\det (B-z\text{\tt Id})}{\det (B-z\text{\tt Id})|_{V_0}}$$
is also  a polynomial with integer coefficients. As shown above, all its roots are of unit modulus. By Kronecker's theorem (\cite{Kronecker}), all these roots are roots of unity.
\ \ \CheckedBox

 \demo{Proof of the displacement lemma}
 \
 
 Let $\{\g_\frak j:\ \frak j\in\mathcal J\}$ be the collection of  eigenvalues of $B|_{V_0^\perp}$ counting multiplicity which are all roots of unity. 
 Let $V_\frak j$ be the corresponding Jordan subspace, then by the above,
 $$\dim\,V_\frak j=2.$$
 We may  extend the period in \Football\ \  as on page \pageref{Football} so that
 $\g_\frak j=1\ \forall\ \frak j\in \mathcal J$.
 \
 
For each $\frak j\in \mathcal J$ let
$(e_j(\frak j):\ j=1,2)$ be the Jordan basis of $V_\frak j$.
\

For $\frak j\in \mathcal J,\ x=x_1e_1(\frak j)+x_2e_2(\frak j)$ and $N\ge 1$, we have that 
$$B^Nx=Nx_1e_1(\frak j)+x_2e_2(\frak j).$$

Thus  for $\Phi:\Bbb Z_Q\to\Bbb R^d\ \&\ 1\le k\le d$,
\begin{align*}
 \s^{(k)}_{K+Ln} & =B^N\s^{(k)}_K\\ &=
 \sum_{\frak j\in\mathcal J}N\<\s^{(k)}_K,e_1(\frak j)\>e_1(\frak j)+\<\s^{(k)}_K,e_2(\frak j)\>e_2(\frak j)\\ &=:
\frak c^{(k)}+N\frak d^{(k)}.
\end{align*}
This proves the displacement lemma.\ \ \Checkedbox

In the sequel, we'll also need the following.
\proclaim{Positivity proposition}\ \ By possibly extending the period in \Football \ as on Page \pageref{Football}, we may ensure that
$B_{s,t}>0\ \forall\ s,t\in S$.\endproclaim
\subsection*{Remark} \ \ Evidently  $E(a(F_{K+1}^{K+L})_{s,t})=\Pi_{F_{K+1}^{K+L}}(0)>0$ iff $B_{s,t}>0$. 
Recall the assumption as in the subsection on subspace decompositions and eigenvalues, that, the parity sequence $(\epsilon_k: k\geq 1)$ is given by:
$$(\e_k:\ k\ge 1)=(\e_1,\dots,\e_K,\overline{\eta_1,\dots,\eta_L}).$$
\demo{Proof}

It follows from (\Pointinghand) as on page \pageref{Pointinghand} that
$$E_{k+1}=
\left(\begin{matrix}
\sum_{j=1}^{m_{k+1}-1}\rho_{\eta_{k}}^{j-1}& \rho_{\eta_{k}}^{m_{k+1}-1}\\
\sum_{j=1}^{m_{k+1}-2}\rho_{\eta_{k}}^{j-1}&  \rho_{\eta_{k}}^{m_{k+1}-2}
\end{matrix}\right).
$$
Choose $1\leq r\leq L$ such that $m_{r}\neq 2$. A direct calculation shows
$$B^2>E_{r+1}E_{r}\geq 
\left(\begin{matrix}
\rho_0+\rho_{\eta_{r}}+\rho_{\eta_{r-1}}& D_1 \\
D_2&D_3
\end{matrix}\right)
$$
where $D_1, D_2, D_3\in M_{Q\times Q}(\N_0)$ are matrices where each row and column has at least one non-zero entry. 

By the parity proposition,  $\eta_k$ and $\eta_{k+1}$ generate the group $\Z_Q$.

Applying this to $k=r-1$, we get that there exists an $N$ such that for all $n\geq N$,
$(\rho_0+\rho_{\eta_{r}}+\rho_{\eta_{r-1}})^n>0$, meaning all of its entries are positive. Thus $B^{N+2}>0$. This proves that $B$ is aperiodic and irreducible and that by extending the period, we can ensure that $B$ is a positive matrix.
\Checkedbox
       
       \subsection*{Asymptotic eventual periodicity $\&$ centering}
   \
       
       \ Let $\a\in\text{\tt QUAD}$ and $\v$ be a step function with rational discontinuities with associated {\tt RAT} sequence $(F_k:\ k\ge 1)$ and
       {\tt ARW} $(X^{(k)}:\ k\ge 1)$.
       
       By the displacement lemma, we may suppose that
  \begin{align*}
    &[n_1,n_2,\dots]=[n_1,\dots,n_K,\overline{m_1,\dots,m_L}];\\ & (\e_k:\ k\ge 1)=(\e_1,\dots,\e_K,\overline{\eta_1,\dots,\eta_L}),\  
    \s_{k+1}=M^{(k+1)}\s_k;\\ &\ (M^{(k)}:\ k\ge 1)=(M^{(1)},\dots,M^{(K)},\overline{E_1,\dots,E_L})\ \&\ \s_{K+Ln}=\frak e+n\frak d.
  \end{align*}

  Next, 
  we examine the  asymptotic, distributional periodicity of the {\tt RAT} sequence and, in particular, that of  the elementary random vector sequence: 
  $$(\frak x^{(k)})=((\frak r^{(k)}_i,\ \frak e^{(k)}_i,\ \frak u^{(k)}(i):\ i=0,1))$$ 
  as on page \pageref{elpres}.
    \proclaim{Elementary  periodic approximation  lemma}
  \ \ 
  
  There are constants $\l,\ M>1$ and, for each $1\le r\le L$ there is a random vector
    $$\frak X^{(r)}:=(\frak R^{(r)}_i,\ \frak E^{(r)}_i,\ \frak U^{(r)}(i):\ i=0,1)\in  RV((\{0,1\}\x\Bbb Z_Q\x\Bbb N_0)^2)$$
    so that

    \begin{align*}&\text{\tt dist}\,(\frak e^{(K+Ln+r)}_i,\ \frak u^{(K+Ln+r)}(i):\ i=0,1\|\frak r^{(K+Ln+r)}_0,\frak r^{(K+Ln+r)}_1)=
  \\ &\ \ \ \ \ \ \ \ \ \ \ \ \ \ \ \ \ \ \ \ \ \ \ \ \ \ \   \text{\tt dist}\,(\frak E^{(r)}_i,\ \frak U^{(r)}(i):\ i=0,1\|\frak R^{(r)}_0,\frak R^{(r)}_1)\text{ and }\\ &
     |P(\frak x^{(K+Ln+r)}=Z)-P(\frak X^{(r)}=Z)|\le \frac{M}{\l^n}\ \ \forall\ n\ge 1,\ Z\in(\{0,1\}\x\Bbb Z_Q\x\Bbb N_0)^2.
    \end{align*} \endproclaim

      \demo{Proof}\ \ We have that
      $$\u\ell_{K+Ln+r}=B(m_r)B(m_{r-1})\cdots B(m_1)C^n\u\ell_K$$
      where
      $$B(m):=\begin{pmatrix} m-1 & 1\\
m-2 & 1
\end{pmatrix}\ \&\ C:=B(m_L)B(m_{L-1})\cdots B(m_1).$$
    Now $\det C=\prod_{r=1}^L\det B(m_r)=1$ and each $C_{i,j}\in\Bbb N$, so $C$  is hyperbolic, with eigenvalues  $\l >1$ and $\frac1\l$.
    \
    
    Moreover, there exists $ c_r(i)\ \ (i=0,1\ \&\ 0\le r\le L)$ with $c_L=\l c_0$ so that
  $$\ell_{K+Ln+r}(i)=c_r(i)\l^n+O(\tfrac1{\l^n}),$$
  whence
   \begin{align*}c_{r+1}(i) & =\frac{\ell_{K+Ln+r+1}(i)}{\l^n}+O(\tfrac1{\l^n})\\ &=
   \frac1{\l^n}[(m_{r+1}-i-1)\ell_{K+Ln+r}(0)+\ell_{K+Ln+r}(1)]+O(\tfrac1{\l^n})\\ &=
    (m_{r+1}-i-1)c_r(0)+c_r(1)+O(\tfrac1{\l^n}).
   \end{align*}
Define random variables $\frak R^{(r)}_i\ \ \ (i=0,1,\ 1\le r\le L)$ by
  $$P(\frak R^{(r+1)}_i=0)=\tfrac{(m_{r+1}-i-1)c_r(0)}{c_{r+1}(i)}\ \&\ P(\frak R^{(r+1)}_i=1)=\tfrac{c_r(1)}{c_{r+1}(i)}=1-P(\frak R^{(r+1)}_i=0).$$
  It follows that for $i,j=0,1\ \&\ 1\le r\le L$,
   $$P(\frak r^{(K+Ln+r)}_i=j)=P(\frak R^{(r)}_i=j)+O(\tfrac1{\l^n}).$$
   Next, we observe that for $n\ge 1,\ 1\le r\le L,\ j=0,1$, the distribution of 
   $\frak e^{(K+Ln+r)}_i$ given $\frak r_i^{(K+Ln+r)}$ does not depend on $n\ge 1$ and define:
   $$P([\frak E^{(r+1)}_i=\D]\|[\frak R^{(r+1)}_i=0])=\frac{\#\frak m(\eta_r,\D)\cap [1,m_{r+1}-i-1]}{m_{r+1}-i-1}\ \ \ (\D\in\Bbb Z_Q)$$ and
   $$P([\frak E^{(r+1)}_i=(m_{r+1}-i-1)\eta_r]\|[\frak R^{(r+1)}_i=1])=1.$$
   Analogously, $\frak u^{(L+Ln+r+1)}(i)$ has a conditional distribution independent of $n$ and we define
$$\frak U^{(r+1)}(i):=\ \ \begin{cases} & \text{uniform on $\frak m(\eta_r,\frak E^{(r+1)}_i)\cap [1,m_{r+1}-i-1]$ if $\frak R_i^{(r+1)}=0$}\\ & 
 m_{r+1}-i\ \ \ \text{if}\ \ \ \frak R_i^{(r+1)}=1.\end{cases}$$
 The random vectors $\frak X^{(r)}\in\text{\tt RV}\,((\{0,1\}\x\Bbb Z_Q\x\Bbb N_0)^2)\ \ \ 1\le r\le L$ where
 \begin{align*}
   \tag{\dsagricultural}\label{dsagricultural}\frak X^{(r)}:=(\frak R^{(r)}_i,\ \frak E^{(r)}_i,\ \frak U^{(r)}(i):\ i=0,1)
 \end{align*}

are as advertised by construction.\ \ \ \Checkedbox

\
 
\proclaim{{\tt RAT} periodic approximation lemma}
 \
 
 There are random variables $a\in\text{\tt RV}\,(M_{S\x S}(\Bbb Z)),\ \frak v,\ \frak w\in \text{\tt RV}\,((\Bbb R^d)^S)$ so that if
       $$H^{(n)}(x)=ax+\frak v+n\frak w\text{ for }x\in (\R^d)^S ,$$ 
       then $\exists\ M>0$ so that $\forall\ n\ge 1,\ f\in M_{S\x S}(\Bbb Z)\x (\Bbb R^d)^S$,        
   \begin{align*}&\tag{\Bicycle}|P(\widetilde{F}_{n}=f)-P(H^{(n)}=f)|\le \frac{M}{\l^n}\\ &\tag{\dsmedical}
   P(\widetilde{F}_{n}=f)>0\ \Leftrightarrow\ P(H^{(n)}=f)>0
   \end{align*}
where
  $$\widetilde{F}_{n}:=F_{K+Ln+1}^{K+Ln+L}.$$
      
 \endproclaim\demo{Proof}
 \
   
   Let $\frak X^{(r)}\ \ (1\le r\le L)$ be independent, each distributed as in (\dsagricultural).

   Define 
   $$\frak l^{(r+1)}_{(i,\e)}:=(\frak R^{(r+1)}_i,\e+\frak E^{(r+1)}_i),$$
  then,  since
   $$\mathcal L^{(K+Ln+r+1)}_{(i,\e)}=(\frak r^{(K+Ln+r+1)}_i,\e+\frak e^{(K+Ln+r+1)}_i),$$
  we have by the elementary periodic approximation lemma,
   $$\sup_{s,t\in S}|P(\frak l^{(r+1)}_s=t)-P(\mathcal L^{(K+Ln+r+1)}_s=t)|=O(\tfrac1{\l^n}).$$

 To study the random variables $W_s^{(K+Ln+r)}$, we'll need formulae for the cumulative displacements.
 \
 
 Using the displacement lemma, for $1\le r\le L$,
  
  \begin{align*} s_{K+Ln+r}(K,\e)&=\sum_{{\nu}=1}^K\s_{K+Ln+r}(0,\e+({\nu}-1)\e_{K+Ln+r})\\ &=\sum_{{\nu}=1}^K(\frak c_r+n\frak d_r)(0,\e+({\nu}-1)\eta_r)
  \\ &=\frak C_r(K,\e)+n\frak D_r(K,\e)
  \end{align*}
  where
  \begin{align*}
\frak C_r(K,\e)&:=\sum_{{\nu}=1}^K\frak c_r(0,\e+({\nu}-1)\eta_r)\\
\frak D_r(K,\e)&:=\sum_{{\nu}=1}^K\frak d_r(0,\e+({\nu}-1)\eta_r).
\end{align*}•

 It follows that

   \begin{align*}\dist(W_{(i,\e)}&^{(K+Ln+r+1)}\|\frak r_i^{(K+Ln+r+1)})\\ &=\dist(s_{K+Ln+r+1}(\frak u^{(K+Ln+r+1)}(i)-1,\e)\|\frak r_i^{(K+Ln+r+1)})\\ &=
   \dist(s_{K+Ln+r+1}(\frak U^{(r+1)}(i)-1,\e)\|\frak R_i^{(r+1)})\\ &=
   \dist(\frak C_{r+1}(\frak U^{(r+1)}(i)-1,\e)+n\frak D_{r+1}(\frak U^{(r+1)}(i)-1,\e)\|\frak R_i^{(r+1)}).
   \end{align*}
   Now let $G^{(n)}_r\ \ (1\le r\le L,\ n\ge 1)$ be the {\tt RAT}s defined by
   $$G^{(n)}_r(x)_{(i,\e)}:=x_{\frak l^{(r)}_{(i,\e)}}+\frak C_{r}(\frak U^{(r)}(i)-1,\e)+n\frak D_{r}(\frak U^{(r)}(i)-1,\e),$$
for all $x\in (\R^d)^S$ then there is a constant $M>0$ so that  $\forall\ f\in M_{S\x S}(\Bbb Z)\x(\Bbb R^d)^S$,
   \begin{align*}
     \tag{\ddag}|P(F_{K+Ln+r}=f)-P(G^{(n)}_r=f)|\le \frac{M}{\l^n}.
    \end{align*}

      Finally, let  
      $$H^{(n)}:=G^{(n)}_{L}\circ G^{(n)}_{L-1}\cdots \circ G^{(n)}_2\circ G^{(n)}_1.$$
   This has the form
       $$H^{(n)}(x)=ax+\frak v+n\frak w$$ 
     for all $x\in (\R^d)^S$ where $a\in\text{\tt RV}\,(M_{S\x S}(\Bbb Z)),\ \frak v,\ \frak w\in \text{\tt RV}\,((\Bbb R^d)^S)$. 
    
       It follows from (\ddag)  that $H^{(n)}$ satisfies (\Bicycle) and (\dsmedical). \ \ \Checkedbox
       
       \subsection*{Coupling}
       \
       
       It follows that there exists a probability space $(\Om,\mathcal A,P)$ on which the independent random vectors $(H^{(n)},\widetilde{F}_n)$ (  $n\ge 1$) can be defined so that
       $$P(H^{(n)}\ne\widetilde{F}_n)\le\frac{M}{\l^n}.$$  \
      Consider the {\tt ARW}
$$Y_J^{(n)}:=H_{J+1}^{n}(X^{(K+LJ)})\ \ \ \ (n>J).$$
\proclaim{{\tt ARW} periodic approximation lemma}

       There is a constant $M>1$ so that for all $\ n>J$,
       \begin{align*}
     \tag{\dscommercial}|P(Y_J^{(n)}\ne X^{(K+Ln)})|\le \frac{M}{\l^J};
    \end{align*}
    and
    \begin{align*}
     \tag{\dsheraldical}\label{dsheraldical}\sup_{n>J}|E(Y_J^{(n)\nu})-E( X^{(K+Ln)\nu})|\xrightarrow[J\to\infty]{}\ 0\ \forall\ \nu\ge 1.
    \end{align*}\endproclaim\demo{Proof of (\dscommercial)}\ \ 
        \begin{align*}
     P(Y_J^{(n)}\ne X^{(K+Ln)})&\le  P(\widetilde{F}_J^n\ne H_J^n) \\ &\le 
     \sum_{{j}=J}^nP(H^{(t)}\ne \widetilde{F}_{j})\\ &\le \sum_{{j}=J}^n\frac{M}{\l^{j}}\\ &=O(\frac1{\l^J}).\ \ \ \CheckedBox
    \end{align*}
    \demo{Proof of (\dsheraldical)}
    For fixed $\nu\ge 1$ and a measurable function $g:\Om\to\Bbb R^S$, for which $|g|^\nu$ is integrable, let
    $$\|g\|_\nu:=E(|g|^\nu)^{\frac1\nu},$$ then $\|\cdot\|_\nu$ is a norm. 
    \
    
    Next,  it follows from the {\tt RAT} periodic approximation lemma that
    $$\|b(\widetilde{F}_n)\|_\infty,\ \|b(H^{(n)})\|_\infty =O(n),$$
        whence
    $$\|Y_J^{(n+1)}\|_\infty,\ \|X^{(K+L(n+1))}\|_\infty=O(n^2).$$ 
    Thus, for some $M'>0$,
\begin{align*} \|Y_J^{(n+1)}&\ -\ X^{(K+L(n+1))}\|_\nu=\|H^{(n+1)}(Y_J^{(n)})\ -\ \widetilde{F}_{n+1}(X^{(K+Ln)})\|_\nu
\\ &\le \|b(H^{(n+1)})\ -\ b(\widetilde{F}_{n+1})\|_\nu+\|(a(H^{(n+1)})\ -\ a(\widetilde{F}_{n+1}))X^{(K+Ln)}\|_\nu+\\ &
\ \ \ \ \ \ \ \ \ \ \ +
\|a(H^{(n+1)})(Y_J^{(n)}\ -\ X^{(K+Ln)})\|_\nu\\ &\le 
\|Y_J^{(n)}\ -\ X^{(K+Ln)}\|_\nu+\frac{M'n^2}{\l^{\frac{n}\nu}}.\end{align*}
Thus possibly increasing $M$,
$$\|Y_J^{(n)}\ -\ X^{(K+Ln))}\|_\nu\le \sum_{j\ge J}\frac{Mj^2}{\l^{\frac{j}\nu}}\xrightarrow[J\to\infty]{}0$$
and (\dsheraldical) follows.\ \ \Checkedbox

\proclaim{Corollary} \ There are constant vectors $\mu,\ \xi,\ \xi_J\in(\Bbb R^d)^S \ \ (J\ge 1)$ and $0<\rho<1$ so that 
      \begin{align*} &
     \  E(Y_J^{({n})})={n}\mu +\xi_J+O(\rho^{n})\ \forall\ J\ge 1,\  \xi_J\xrightarrow[J\to\infty]{}\ \xi,\\ &\&\ \ 
      E(X^{(K+L{n})})={n}\mu +\xi+O(\rho^{n}).     
      \end{align*} \endproclaim

 \demo{Proof}\ \   We have 
  $$Y_J^{({n}+1)}=H^{({n}+1)}(Y_J^{({n})})$$ where  $$H^{(n)}(x)=a^{(n)}x+\frak v^{(n)}+n\frak w^{(n)}$$ and
  $(a^{(n)},\frak v^{(n)},\frak w^{(n)}:\ n\ge 1)$ are independent and  identically distributed.

\

It follows as in \cite{ABN2016} that
  \begin{align*}\tag{\dschemical}\label{dschemical}E(Y_J^{({n})})=E(a)^{n}E(X^{(K+LJ)})+\sum_{k =1}^{n}E(a)^{{n}-k}E(\frak v+k\frak w)
     \end{align*}
     By the positivity proposition, by possibly extending the period in \Football \ as on page \pageref{Football}, we may ensure that 
$E(a(H^{(n)}))=\Pi_{ H^{(n)}}(0)$ is an aperiodic stochastic matrix whence $1$ is a simple, dominant eigenvalue with eigenvector $\mathbb{1}\in\Bbb C^S$.

\

Suppose that $\pi\in\Bbb R_+^S$ satisfies $\<\pi,\mathbb{1}\>=1\ \&\ E(a)^*\pi=\pi$ (where $A*$ is the transpose of the  matrix $A$).

\

Let $N:\Bbb C^S\to\Bbb C\cdot\mathbb{1},\ N(x):=\<\pi,x\>\mathbb{1}$, then 
$$E(a)^nx=N(x)\mathbb{1}+R^nx$$
where $Rx:=E(a)(x-N(x))$ and $\exists\ 0<\rho<1$ so that $\|R^n\|=O(\rho^n)$.

We claim next that $\exists\ \mu,\ \xi_J\ \in(\Bbb R^d)^S$ so that
\begin{align*}\tag{\Radioactivity} E(Y_J^{({n})})=\frac{{n}({n}+1)}2N(E(\frak w))+ {n}\mu +\xi_J+O(t\rho^{n}). 
\end{align*}
\demo{Proof of (\Radioactivity)}
\

By (\dschemical) as on page \pageref{dschemical},
\begin{align*}E(Y_J^{({n})})&=E(a)^{n}E(X^{(K+LJ)})+\sum_{k=1}^{n}E(a)^{{n}-k}E(\frak v+k\frak w)\\ &=
N(E(X^{(K+LJ)}))+({n}-1)N(E(\frak v))+\frac{{n}({n}-1)}2N(E(\frak w))+\mathcal E({n})
     \end{align*}
     where
     
\begin{align*}\mathcal E({n})&:=R^{n}E(X^{(K+LJ)})+\sum_{k=1}^{n}R^{{n}-k}E(\frak v+k\frak w)\\ &
=\sum_{{n}=1}^{n}R^{{n}-k}E(\frak v+k\frak w)+O(\rho^{n})\\ &=\sum_{{n}=0}^{k-1}({n}-k)R^k(\frak w)+\sum_{k=0}^{{n}-1}R^k E(\frak v)+O(\rho^{n})\\ &=
t\sum_{k=0}^\infty R^k(\frak w)-\sum_{k=0}^\infty k R^k E(\frak w)+\sum_{k=0}^\infty R^k E(\frak v)+O({n}\rho^{n}).\ \ \CheckedBox\ \text{(\Radioactivity)}
\end{align*}

 To obtain the expansion for $E(Y_J^{({n})})$ from (\Radioactivity) (with enlarged $\rho$), it suffices to show 
 that $N(E(\frak w))=0$.
 \
 
 This will follow from the Denjoy-Koksma estimate.
 \
 
 By (\dsheraldical) as on page \pageref{dsheraldical}, we have
 $$|E(X^{(K+L{n})})-E(Y_J^{({n})})|=O(1).$$
Thus, if 
$N(E(\frak w))\ne 0$, then by (\Radioactivity),  $|E(X^{(K+L{n})})|\asymp {n}^2$ contradicting the Denjoy-Koksma estimate that $|E(X^{(K+L{n})})|=O({n})$. 
The expansion for $E(X^{(K+L{n})})$  follows.\ \ \Checkedbox
  
  \subsection*{Centering}
  \
  \label{page: fancy F}
  As in \cite{ABN2016}, set $(\widehat{X}^{(n)}:\ n\ge 1)$ be the centered {\tt ARW} defined by 
   $$\widehat{X}^{(n)}:={X}^{(n)}-E(X^{(n)})$$ and let $(\F_n:\ n\ge 1)$ be the independent {\tt RAT} sequence so that
   $$\widehat{X}^{(K+Ln)}=\F_1^n(\widehat{X}^{(K)}).$$
 
\

\proclaim{{\tt ARW} centering lemma}
  \
  
  There is a centered,  independent, identically distributed {\tt RAT} sequence $(\mathcal H_n:\ n\ge 1)$ and $0<\rho<r<1$ so that if 
 for $J\ge 1$,  
 $(Z_J^{({n})}:\ {n} > J)$ is defined by
 $$Z_J^{({n})}:=\mathcal H_{J+1}^{n}(\widehat{X}^{(K+LJ)}),$$ then
   \begin{align*}&\tag{i}\sup_{n>J}|E(Z_J^{(n)\nu})-E(\widehat{X}^{(K+Ln)\nu})|\xrightarrow[J\to\infty]{}\ 0\ \forall\ \nu\ge 1;\\ &
    \tag{ii} P(\exists\ {n}\ge J,\ |Z_J^{({n})}- \widehat{X}^{(K+L{n})}|\ge r^{n})=O(\rho^J)\ \text{as}\ J\to\infty.
   \end{align*}

       \endproclaim\demo{Proof}
       \
       
      Define  $$\widehat{Y}_J^{({n})}:=Y_J^{({n})}-E(Y_J^{({n})})=:Y_J^{({n})}-c_{n}.$$ 
  As in \cite{ABN2016}, $(\widehat{Y}_J^{({n})}:\ {n}\ge 1)$ is given by the centered {\tt RAT} sequence $(\mathcal G_n:\ n\ge 1)$ where
  $a(\mathcal G_n)=a^{(n)}$ and
  \begin{align*}b(\mathcal G_{n+1})&= a^{(n+1)}c_n-c_{n+1}+\frak v^{(n+1)}+n\frak w^{(n+1)} \\ &=
  (a^{(n+1)}-I)\xi_J+\frak v^{(n+1)}-\mu+ n[(a^{(n+1)}-I)\mu+\frak w^{(n+1)}]+O(\rho^n)\\ &=:
  \frak v{'}^{(n+1)}+n\frak w{'}^{(n+1)}+O(\rho^n)
  \end{align*}
where $(a^{(n)},\frak v^{(n)},\frak w^{(n)})$ are independent, identically distributed random variables and $I$ is the identity matrix.
\

By the remark after the positivity proposition,  $E(a)$ is irreducible and aperiodic. 
\

Thus, by the variance lemma in \cite{ABN2016}, for each $s\in S$,
\begin{align*}E((\widehat{Y}_J^{({n})2})_s)&\asymp\ \sum_{k=1}^{n}E(b(\mathcal G_k)_s^2)\\ &=
\sum_{k=1}^{n}E((\frak v'_s+k\frak w'_s+O(\rho^k))^2)\\ &\sim
{n}E(\frak v_s^{\prime 2})+{n}^2E(\frak w'_s\frak v'_s)+\frac{{n}^3}3E(\frak w_s^{\prime 2}).
   \end{align*}
   
   \
   
By (\dsheraldical) as on page \pageref{dsheraldical},
   $$|E((\widehat{Y}_J^{({n})2})_s)-E((\widehat{X}^{(K+L{n})2})_s)|=O(1)$$
   whence, by the Denjoy-Koksma estimate, $E((\widehat{Y}_J^{({n})2})_s)\ll {n}^2$ and  $\frak w{'}^{(n+1)}\equiv 0$.
\

Thus 
   \begin{align*}b(\mathcal G_{n+1})= 
  (a^{(n+1)}-I)\xi_J+\frak v^{(n+1)}-\mu+O(\rho^n).
  \end{align*}
  Accordingly, define $\mathcal H_n$ by
  $$a(\mathcal H_n):=a^{(n)}\ \&\ b(\mathcal H_n):=(a^{(n)}-I)\xi_J+\frak v^{(n)}-\mu.$$
   The lemma follows.\ \ \Checkedbox 
   \section*{\S5  Spectral theory and theorem 2}
   \

  By the Perron-Frobenius theorem, $1$ is a simple, dominant eigenvalue of $\Pi_{\mathcal H}(0)$ (where $\mathcal H:=\mathcal H_1$  and 
  $\Pi_{\mathcal H}$ is the {\tt RAT-CF} as defined by (\dsliterary) on page \pageref{lit}) 
  with  right
 eigenvector $\mathbb{1}\in\Bbb R_+^S$ and left eigenvector $\pi\in\Bbb R_+^S$ satisfying $\<\pi,\mathbb{1}\>=1$.
 \
 
 By the implicit function theorem $\exists\ r=r_{\mathcal H}>0$ and smooth functions 
 $$\l:(-r,r)^d\to\Bbb C,\ v:(-r,r)^d\to\Bbb C^S,\ \pi:(-r,r)^d\to\Bbb C^S$$ so that 
\bul  $\<\pi(0),v(\th)\>=\<\pi(\th),v(\th)\>=1$;
  
  \bul $\l(0)=1,\ v(0)=\mathbb{1}\ \&\ \pi(0)=\pi$;
   \bul for each $1\le k\le\infty,\ \th\in (-r,r)^d$,\  $\l(\th)$ is a simple, dominant eigenvalue of $\Pi_{\mathcal H}(\th)$ with eigenvector $v(\th)$ and left eigenvector
   $\pi(\th)$.

      As in \cite{HH}, consider the {\it principal projections} $N(\th):\Bbb C^S\to\Bbb C^S$ defined by
   $$N(\th)x:=\<\pi(\th),x\>v(\th)$$ then possibly reducing $r_\mathcal H>0$, we ensure $\exists\ 0<\rho<1$ so that
   $$\Pi_{\mathcal H}(\th)^n-\l(\th)^nN(\th)=R(\th)^n=O(\rho^n)\ \text{uniformly in}\ |\th|\le r_\mathcal H$$
  where $R(\th):=\Pi_{\mathcal H}(\th)(I-\l(\th))N(\th)$.
 
 \
 
 The proofs of our limit theorems in the sequel  use the following lemma.
     
\proclaim{Lemma:\ {\tt Taylor expansion of the eigenvalue}}\ \ Under the assumptions of Theroem 2,
\begin{align*}\l(\th)=1-\<D\th,\th\>+o(\|\th\|^2)
\end{align*}

as $\th\to\ 0$ where $D\in M_{d\x d}(\Bbb C)$ is  positive definite.\endproclaim
\demo{Proof}\ \   We have
   $$\l(\th)=1+\<\nabla\l(0),\th\>+\<d^2\l(0)\th,\th\>+o(\|\th\|^2)$$
   where $d^2\l(\th)$ is the matrix of second partial derivatives: $$d^2\l(\th)_{h,j}:=\frac{\bdy^2\l}{\bdy\th_h\bdy\th_j}(\th).$$
   and we must show that 
   $\nabla\l(0)=0$ and that $D:=-d^2\l(0)$ is positive definite.

   \

Fix  $\s\in\Bbb R^d,\ \|\s\|=1$ and write, for differentiable $f:\Bbb R^d\to\Bbb C$,
$$D_\s^k f(\th):=\tfrac{d^k}{dt^k}f(\th+t\s)|_{t=0},$$ then
\begin{align*}D_\s f(\th)=\<\s,\nabla f(\th)\>\ \&\ D_\s^2 f(\th)=\<d^2f(\th)\s,\s\>. 
\end{align*}
Accordingly, it suffices to show that for each $\s\in\Bbb R^d,\ \|\s\|=1,$
   \sms (i) $D_\s\l(0)=\tfrac{d}{dt}\l(t\s)|_{t=0}=0$ and (ii)  $D_\s^2\l(0)=\tfrac{d^2}{dt^2}\l(t\s)|_{t=0}<0$.
\Par1 $D_\s(\Pi_\mathcal H)(0)\mathbb{1}=0$.
\demo{Proof of \P1}\ \ For fixed $s\in S$: 

\begin{align*}
    (D_\s(\Pi_\mathcal H)(0)\mathbb{1})_s &=\frak i\sum_{t\in S}P(\mathcal L_s(\mathcal H)=t)E(\<\s,b_s(\mathcal H)\>\|\mathcal L_s(\mathcal H)=t)\\ &=\frak iE(\<\s,b_s(\mathcal H)\>)\\ &=
      0\ \ \because\ \mathcal H\ \ \text{is centered.\ \ \Checkedbox \ \ \P1}
   \end{align*}
   
   \demo{Proof of {\rm (i)}}

Since $\<\pi(0),v(\th)\>\equiv 1$, we have that $D_\s v(\th)\perp\pi(0)$. Also
\begin{align*}&D_\s(\Pi_\mathcal Hv)(\th)=D_\s(\Pi_\mathcal H)(\th)v(\th)+\Pi_\mathcal H(\th)D_\s v(\th),\\ &
D_\s(\l v)(\th)=D_\s\l(\th)v(\th)+\l(\th)D_\s v(\th)
\end{align*}
 whence
  \begin{align*}0&=D_\s (\Pi_\mathcal Hv-\l v)(\th)\\ &=
  D_\s(\Pi_\mathcal H-\l)(\th)v(\th)+(\Pi_\mathcal H-\l)(\th)D_\s v(\th)   
  \end{align*}
  and in particular 
  \begin{align*}0&=\<\pi(0),D_\s (\Pi_\mathcal Hv-\l v)(0)\>\\ &=
  \<\pi(0),D_\s(\Pi_\mathcal H-\l)(0)\mathbb{1}\>+\<\pi(0),(\Pi_\mathcal H(0)-1)D_\s v(0)\>  \\ &=
  \<\pi(0),D_\s(\Pi_\mathcal H-\l)(0)\mathbb{1}\>\ \because\ \Pi_\mathcal H(0)^*\pi(0)=\pi(0).
  \end{align*}
  Thus
   \begin{align*}
      D_\s\l(0)=\<\pi(0),D_\s(\Pi_\mathcal H)(0)\mathbb{1}\>=0\ \ \text{by \ \P1.\ \ \Checkedbox\ \ (i)}
   \end{align*}
   
   \Par2 $D_\s v(0)=0$.
   \demo{Proof}\ \ Differentiating $\Pi_\mathcal H(\th)v(\th)=\l(\th)v(\th)$ at $0$:
  $$D_\s(\Pi_\mathcal H)\mathbb{1}+\Pi_\mathcal H(0)D_\s v(0)=\ D_\s\l(0)\mathbb{1}+D_\s v(0).$$
  By \P1 $\&$ (i),
  $$\Pi_\mathcal H(0)D_\s v(0)=D_\s v(0).$$
  Thus $D_\s v(0)=c\mathbb{1}$ for some $c\in\Bbb C$. But $D_\s v(0)\perp\pi(0)$, and so
  $$\overline{c}=\<\pi(0),D_\s v(0)\>=0. \ \ \  \CheckedBox\ \ \P2$$
  \Par3 $D_\s^2(\l)(0)=\<\pi(0),D_\s^2(\Pi_\mathcal H)(0)\mathbb{1}\>$.
  \demo{Proof}\ \ Differentiating $\Pi_\mathcal H(\th)v(\th)=\l(\th)v(\th)$ twice at $0$ in direction $\s$:
\begin{align*}
 D_\s^2(\Pi_\mathcal H)(0)&\mathbb{1}+2D_\s(\Pi_\mathcal H)(0)D_\s v(0)+\Pi_\mathcal H(0)D_\s^2 v(0)\\& =D_\s^2(\l)(0)\mathbb{1}+2D_\s(\l)(0)D_\s v(0)+D_\s^2 v(0).
\end{align*}
  By (i) $\&$ \P2
    $$D_\s^2(\Pi_\mathcal H)(0)\mathbb{1}+\Pi_\mathcal H(0)D_\s^2 v(0)=D_\s^2(\l)(0)\mathbb{1}+D_\s^2 v(0)$$
    and in particular
 \begin{align*}D_\s^2(\l)(0)&=\<\pi(0),D_\s^2(\l)(0)\mathbb{1}\>\\ &=
 \<\pi(0),D_\s^2(\Pi_\mathcal H)(0)\mathbb{1}\>+\<\pi(0),\Pi_\mathcal H(0)D_\s^2 v(0)\>-\<\pi(0),D_\s^2 v(0)\>\\ &=
 \<\pi(0),D_\s^2(\Pi_\mathcal H)(0)\mathbb{1}\>\ \ \because\ \Pi_\mathcal H(0)^*\pi(0)=\pi(0).\ \ \ \CheckedBox\ \ \P3
  \end{align*}
  \demo{Proof of {\rm (ii)}}
\ \ For fixed $s\in S$:
\begin{align*}
    (D_\s^2(\Pi_\mathcal H)(0)\mathbb{1})_s &=-\sum_{t\in S}P(\mathcal L_s(\mathcal H)=t)E(\<\s,b_s(\mathcal H)\>^2\|\mathcal L_s(\mathcal H)=t)\\ &=-E(\<\s,b_s(\mathcal H)\>^2)
   \end{align*}
   whence by \P3,
   \begin{align*}D_\s^2(\l)(0)& =\<\pi(0),D_\s^2(\Pi_\mathcal H)(0)\mathbb{1}\>\\ &=-\sum_{s\in S}\pi_s(0)E(\<\s,b_s(\mathcal H)\>^2)\ \le 0    
   \end{align*}
with equality iff $\<\s,b_s(\mathcal H)\>=0\ \forall\ s\in S$.

\

Next recall that $Z^{(n)}:=\mathcal H_1^n(0)=X^{(K+Ln)}-n\mu +O(1)$.

\

If $\<\s,b_s(\mathcal H)\>=0\ \forall\ s\in S$, then, taking $s=(0,0)$ we have
$$\sup_n|\<\s,Z_s^{(n)}\>|<\infty\ \Lra \ \sup_n|\<\s,X^{(K+Ln)}_s-n\mu_s \>|=:W<\infty$$
whence

$$|\<\s,\v_j(0)-n\mu_s \>|\le W\ \forall\ n\ge 1,\ 1\le j\le \ell_{K+Ln}(0).$$
It follows from this that $\mu_s=0$ and  that $\<\s,\v\>$ is  a coboundary. 
\

In view of the assumption that $\text{\tt span}_\Bbb R\,\v(\Bbb T)$ has full dimension, this contradicts theorem 1 and completes the proof of the lemma.\ \ \Checkedbox\ \ (ii)\ \qed
   \demo{Proof of theorem 2} \ \ It suffices to prove that for fixed $s\in S,\ \th\in\Bbb R^d$ 
\begin{align*}\tag{\dsrailways}E(\exp[2\pi\frak i \<\th,&\tfrac{X^{(K+Ln)}_s-n\mu_s}{\sqrt n}])\xrightarrow[n\to\infty]{}
\exp[-\tfrac{\<\th,D\th\>}2]. 
\end{align*}

By asymptotic, eventual periodicity $\exists\ \rho\in (0,1)$ so 
that for any fixed $J,\ s\in S$ and all $n\ge 1$,
\begin{align*}E(\exp[2\pi\frak i \<\th,&\tfrac{X^{(K+L(J+n))}_s-(J+n)\mu_s}{\sqrt n}\>])=E(\exp[2\pi\frak i \<\th,\tfrac{(Z{_J}^{(n)})_s}{\sqrt n}\>])+O_n(\rho^J)\\ &=
        \l(\tfrac{\th}{\sqrt n})^nE(\exp[2\pi\frak i \<\th,\tfrac{X^{(K+LJ)}_s-J\mu_s}{\sqrt n}\>])+O_n(\rho^J)
     \end{align*}
     where $|O_n(\rho^J)|\leq M \rho	^J$ for a constant $M$ independent of $n$.
     By the Taylor expansion of the eigenvalue,
   \begin{align*}\l(\tfrac{\th}{\sqrt n})^nE(\exp[2\pi\frak i \<\th,\tfrac{X^{(K+LJ)}_s-J\mu_s}{\sqrt n}\>])&
   \xrightarrow[n\to\infty]{}\ \exp[-\tfrac{\<\th,D\th\>}2]\ \forall\ J\ge 1.    
   \end{align*}

     To deduce (\dsrailways)  from this, let $\varepsilon>0$ and choose $J=J_\varepsilon\ge 1$ so that
     $$|E(\exp[2\pi\frak i \<\th,\tfrac{X^{(K+L(J+n))}_s-(J+n)\mu_s}{\sqrt n}\>])-
        \l(\tfrac{\th}{\sqrt n})^nE(\exp[2\pi\frak i \<\th,\tfrac{X^{(K+LJ)}_s-J\mu_s}{\sqrt n}\>])|<\frac{\varepsilon}2\ \forall\ n\ge J$$
     and then choose $N=N_{J,\varepsilon}>J$ so that
    $$|\l(\tfrac{\th}{\sqrt n})^nE(\exp[2\pi\frak i \<\th,\tfrac{X^{(K+LJ)}_s-J\mu_s}{\sqrt n}\>])-\exp[-\tfrac{\<\th,D\th\>}2]|<\frac{\varepsilon}2
    \ \forall\ n>N.$$
    This implies (\dsrailways). \ \ \qed
   
   \section*{\S6   The {\tt WRLLT} and theorem 3}

   \
   
   \subsection*{Visits to zero and {\tt RAT}{\tt s}}
\ \ \ Recall that we assume $Q\ge 1,\ \a\in[0,1)\setminus\Bbb Q$, with $\{Q\a\}=[n_1,n_2,\dots]$.
 \
 
 Let $\Phi:\Bbb Z_Q\to\Bbb Z^d$ satisfy $\sum_{k\in\Bbb Z_Q}\Phi(k)=0\ \&\ \<\Phi(\Bbb Z_Q)\>=\Bbb Z^d$ 
and define $\v:\Bbb T\to\Bbb Z^d$ by
$$\v(x):=\Phi(\lfl Qx\rfl)$$ and $T=T_{\a,\Phi}:\Bbb T\x\Bbb Z^d\to\Bbb T\x\Bbb Z^d$ by
 $$T(x,z):=(\{x+\a\},z+\v(x)).$$
 \
 Let $ \Psi_N(x):=\#\{1\le n\le N:\ \v_n(x)=0\}$.
  \proclaim{Visit lemma }\ \ \ Let $(X^{(k)}:\ k\ge 1)$ be the associated {\tt ARW}, then
\begin{align*}&\tag{a} \int_0^1\Psi_{\ell_k(0)}(x)dx\  \ge\ \ 
\frac{\ell_k(1)^2}{3\ell_k(0)}\min_{\e\in\Bbb Z_Q}\int_{\mathbb T^d}|E(e^{2\pi\frak i \<\th, X^{(k)}(1,\e)\>})|^2d\th-\frac12
\\ &\tag{b}\|\Psi_{\ell_k(1)}\|_\infty\le
2\ell_k(0)\max_{\e\in\Bbb Z_Q}\int_{\mathbb T^d}|E(e^{2\pi\frak i \<\th, X^{(k)}(0,\e)\>})|d\th.
\end{align*}\endproclaim

\subsection*{Visit sets}
 \
 
 The {\it visit set to $\nu\in\mathbb Z^d$} is
 $$K_\nu:=\{n\ge 1:\ \v_n(0)=\nu\}$$
 and the {\it visit distributions} are the discrete  measures $U_k^{(i)}$ on $\mathbb Z^d$ defined by
 $$U_k^{(i)}(\nu):=\#(K_\nu\cap [1,\ell_k(i)])\ \ \ (k\ge 1,\ i=0,1).$$
 The {\it auxiliary visit sets to $\nu\in\mathbb Z^d$} are
 $$K_k(i,\e,\nu):=\{1\le j\le \ell_k(i):\ \Si_k(i,\e)_j=\nu\}.$$
 and the {\it auxiliary visit distributions} are the discrete  measures $V_k(i,\e)$ on $\mathbb Z^d$ defined by
 $$V_k(i,\e)(\nu):=\#(K_k(i,\e,\nu))\ \ \ (k\ge 1,\ i=0,1).$$
 As above,
 $$K_k(0,0,\nu)=K_\nu\cap [1,\ell_k(0)]\ \&\ \ U_k^{(0)}=V_k(0,0).$$
 \proclaim{Visit sublemma }
 \begin{align*}
 &\tag{4.1}  \int_0^1\Psi_{\ell_k(0)}(x)dx\  \ge\ \ \frac1{3\ell_k(0)}\min_{\e\in\Bbb Z_Q}\sum_{\nu\in\mathbb Z^d}V_k(1,\e)(\nu)^2-\frac12;\\ &\tag{4.2}
 \int_0^1\Psi_{\ell_k(1)}(x)^Ndx\  \le\ \ \frac{2^N}{\ell_k(1)}\max_{\e\in\Bbb Z_Q}\sum_{\nu\in\mathbb Z^d}V_{k}(0,\e)(\nu)^{N+1}\ \ \forall\ N\ge 1.
 \end{align*}
 \endproclaim
\demo{Proof}
Fix $N,\ k\ge 1\ \&\ i=0,1,\ \Psi_{\ell_k(i)}^N$ is a step function on $\Bbb T$, whence Riemann integrable. Using the unique ergodicity of $x\mapsto x+\a$,
 
 \begin{align*}\tag*{\Yinyang}\label{Yinyang}\ell_{k+r}(0)\int_0^1&\Psi_{\ell_k(i)}(x)^Ndx\\ &
 \underset{r\to\infty}{\text{\Large $\sim$}}\sum_{j=1}^{\ell_{k+r}(0)}\Psi_{\ell_k(i)}(\{j\a\})^N\\ &=
 \sum_{j=1}^{\ell_{k+r}(0)}\#(\{1\le m\le \ell_k(i):\ \v_m(\{j\a\})=0\})^N\\ &=
 \sum_{j=1}^{\ell_{k+r}(0)}\#(\{1\le m\le \ell_k(i):\ \v_{m+j}(0)=\v_j(0)\})^N\\ &=
 \sum_{j=1}^{\ell_{k+r}(0)}\#(\{j+1\le m\le \ell_k(i)+j:\ \v_{m}(0)=\v_j(0)\})^N\\ &=
 \sum_{\nu\in\mathbb Z^d}\sum_{j\in [1,\ell_{k+r}(0)]\cap K_\nu}\#\,(K_\nu\cap [j+1,j+\ell_k(i)])^N. \end{align*}
 
By \cite[Theorem 2.1]{AK} and the orbit block transitions, for $r\ge 1,\ \exists\ J=J_{r,k}\ge 1,\ 0=m_1<\dots<m_J=\ell_{k+r}(0)$ and $\eta_1,\dots\eta_{J-1}\in\Bbb Z_Q,\ i_1,\dots,i_{J-1}=0,1$  so that 
$$m_{j+1}-m_j=\ell_k(i_j)\ \forall\ j,\ \ [1,\ell_{k+r}(0)]=\bigcupdot_{j=1}^{J-1}(m_j,m_{j+1}]$$ 
and 
\begin{align*}
 \tag{\dsmilitary}\label{dsmilitary}(\v_{m_j+1}(0),\v_{m_j+2}(0),\dots,\v_{m_{j+1}}(0))=\v_{m_{j}}(0)+\Si_k(i_j,\eta_j).   \end{align*}
Since $m_{j+1}-m_j=\ell_k(i_j)$, it follows that 
 \begin{align*}
  \tag*{\WashCotton}\ell_{k+r}(0)\in [J\ell_k(1),J\ell_k(0)].
 \end{align*}
\

Also by (\dsmilitary), for fixed $\nu\in\mathbb Z^d,$
 \begin{align*}
  K_\nu\cap (m_j,m_{j+1}]=m_j+K_k(i_j,\eta_j,\nu-\v_{m_j}(0)).
 \end{align*}
\demo{Proof of {\rm (4.1)}}\  \ We have for fixed $\nu\in\mathbb Z^d$
 \begin{align*}
\sum&_{h\in [1,\ell_{k+r}(0)]\cap K_\nu}\#\,(K_\nu\cap [h+1,h+\ell_k(0)])\\ &=
\sum_{j=1}^{J-1}\sum_{h\in (m_j,m_{j+1}]\cap K_\nu}\#\,(K_\nu\cap [h+1,h+\ell_k(0)])\\ &\ge
\sum_{j=1}^{J-1}\sum_{h\in (m_j,m_{j+1}]\cap K_\nu}
\#\,(K_\nu\cap [h+1,m_{j+1}])\ \ \because\ \forall\ h,j:\ h+\ell_k(0)\ge h+\ell_k(i_j)\ge m_{j+1}\\ &\ge
\sum_{j=1}^{J-1}\sum_{k=1}^{L_j-1}k\ \ \text{where}\ \  L_j:=\#((m_j,m_{j+1}]\cap K_\nu)=V_k(i_j,\eta_j)(\nu-\v_{m_j}(0))
\\ &=\sum_{j=1}^{J-1}\frak t(L_j)\ \text{where}\ \frak t(x):=\frac{x(x-1)}2\\ &=
\sum_{j=1}^{J-1}\frak t(V_k(i_j,\eta_j)(\nu-\v_{m_j}(0))).
\end{align*}
Thus
\begin{align*}
  \sum_{\nu\in\mathbb Z^d}&\sum_{j\in [1,\ell_{k+r}(0)]\cap K_\nu}\#\,(K_\nu\cap [j+1,j+\ell_k(0)])\\ &\ge
 \sum_{\nu\in\mathbb Z^d}\sum_{j=1}^{J-1}\frak t(V_k(i_j,\eta_j)(\nu-\v_{m_j}(0)))
 \\ &\ge (J-1)\min_{\e\in\Bbb Z_Q}\sum_{\nu\in\mathbb Z^d}\frak t(V_k(1,\e)(\nu))\ \ \text{since $\frak t\ \uparrow$ on $\Bbb N_0$}
 \\ &=
 \frac{J-1}2\min_{\e\in\Bbb Z_Q}\sum_{\nu\in\mathbb Z^d}V_k(1,\e)(\nu)^2-\frac{(J-1)\ell_k(1)}2
\end{align*}
and using \Yinyang\  as on page \pageref{Yinyang} with $N=1\ \&\ i=0$ 
 \begin{align*} \int_0^1&\Psi_{\ell_k(0)}(x)dx \xleftarrow[r\to\infty]{}\frac1{\ell_{k+r}(0)}
  \sum_{j=1}^{\ell_{k+r}(0)}\Psi_{\ell_k(0)}(\{j\a\})\\ &\ge
  \frac{J-1}{2\ell_{k+r}(0)}\min_{\e\in\Bbb Z_Q}\sum_{\nu\in\mathbb Z^d}V_k(1,\e)(\nu)^2-\frac{J\ell_k(1)}{2\ell_{k+r}(0)}
  \\ &\ge \frac1{3\ell_{k}(0)}\min_{\e\in\Bbb Z_Q}\sum_{\nu\in\mathbb Z^d}V_k(1,\e)(\nu)^2-\frac12\ \ \text{by \WashCotton.\ \  \Checkedbox\ (4.1)}.
 \end{align*}

\demo{Proof of {\rm (4.2)}}

\

\ \ Using \Yinyang\ as on page \pageref{Yinyang} with $k,\ N\ge 1$ arbitrary and fixed $\ \&\ i=1$ we have

 \begin{align*}\ell_{k+r}(0)\int_0^1\Psi_{\ell_k(1)}(x)^Ndx&
 \underset{r\to\infty}{\text{\Large $\sim$}}  \sum_{\nu\in\mathbb Z^d}\sum_{j\in [1,\ell_{k+r}(0)]\cap K_\nu}\#\,(K_\nu\cap [j+1,j+\ell_k(1)])^N.\end{align*}
Similar to (\dsmilitary) as on page \pageref{dsmilitary},
\begin{align*}
(\v_{m_j+1}(0),\v_{m_j+2}(0),\dots,\v_{m_{j+1}+\ell_k(1)}(0))=[\v_{m_{j}}(0)\mathbb{1}+\Si_k(i_j,\eta_j)]\odot [\v_{m_{j+1}}(0)\mathbb{1}+\Si_k(1,\D_j)].   \end{align*}
 for some $\D_j\in\Bbb Z_Q$.

 We have 
as before, for fixed $\nu\in\mathbb Z^d$,
\begin{align*}
\sum_{h\in [1,\ell_{k+r}(0)]\cap K_\nu}\#\,(K_\nu\cap [h+1,h+\ell_k(1)])^N&=
\sum_{j=1}^{J-1}\sum_{h\in (m_j,m_{j+1}]\cap K_\nu}\#\,(K_\nu\cap [h+1,h+\ell_k(1)])^N\\ &\le
\sum_{j=1}^{J-1}\sum_{h\in (m_j,m_{j+1}]\cap K_\nu}\#\,(K_\nu\cap [h+1,m_{j+1}+\ell_k(1)])^N.
\end{align*}
Fix $j$. For fixed $h\in (m_j,m_{j+1}]$,
\begin{align*}\#\,(K_\nu\cap [h+1,m_{j+1}+\ell_k(1))])&=\#\,(K_\nu\cap [h+1,m_{j+1}])+\#\,(K_\nu\cap [m_{j+1}+1,m_{j+1}+\ell_k(1))])\\ &=
\#\,(K_\nu\cap [h+1,m_{j+1}])+V_k(1,\D_j)(\nu-\v_{m_{j+1}}(0)). 
\end{align*}
Thus 
\begin{align*}
 &\sum_{h\in (m_j,m_{j+1}]\cap K_\nu}\#\,(K_\nu\cap [h+1,m_{j+1}+\ell_k(1)])^N=
 \\ &=\sum_{h\in (m_j,m_{j+1}]\cap K_\nu}\#\,(K_\nu\cap [h+1,m_{j+1}])+V_k(1,\D_j)(\nu-\v_{m_{j+1}}(0))^N
 \\ &=\sum_{n=0}^N\binom{N}{n}\(\sum_{h\in (m_j,m_{j+1}]\cap K_\nu}\#\,(K_\nu\cap [h+1,m_{j+1}])^n\)V_k(1,\D_j)(\nu-\v_{m_{j+1}}(0))^{N-n}
 \\ &\le\sum_{n=0}^N\binom{N}n V_k(i_j,\eta_j)(\nu-\v_{m_{j}}(0))^{n+1}V_k(1,\D_j)(\nu-\v_{m_{j+1}}(0))^{N-n}
 \\ &\le\sum_{n=0}^N\binom{N}{n} V_k(0,\eta_j)(\nu-\v_{m_{j}}(0))^{n+1}V_k(0,\D_j)(\nu-\v_{m_{j+1}}(0))^{N-n}.
\end{align*}

Using this and H\"older's inequality,
\begin{align*}
 &\sum_{\nu\in\mathbb Z^d}\sum_{j\in [1,\ell_{k+r}(0)]\cap K_\nu}\#\,(K_\nu\cap [j+1,j+\ell_k(1)])^N\le 
 \\ &\le\sum_{j=1}^J\sum_{n=0}^N\binom{N}{n}\sum_{\nu\in\mathbb Z^d} V_k(0,\eta_j)(\nu-\v_{m_{j}}(0))^{n+1}V_k(0,\D_j)(\nu-\v_{m_{j+1}}(0))^{N-n}
  \\ &\le\sum_{j=1}^J\sum_{n=0}^N\binom{N}{n}\(\sum_{\nu\in\mathbb Z^d} V_k(0,\eta_j)(\nu-\v_{m_{j}}(0))^{N+1}\)^{\frac{n+1}{N+1}}
  \(\sum_{\nu\in\mathbb Z^d}V_k(0,\D_j)(\nu-\v_{m_{j+1}}(0))^{N+1}\)^{\frac{N-n}{N+1}}
  \\ &=2^NJ\max_{\e\in\Bbb Z_Q}\sum_{\nu\in\mathbb Z^d}V_k(0,\e)(\nu)^{N+1}
  \end{align*}
              
whence
\begin{align*}\int_0^1\Psi_{\ell_k(1)}(x)^Ndx&
 \xleftarrow[r\to\infty]{}\frac1{\ell_{k+r}(0)}\sum_{\nu\in\mathbb Z^d}\sum_{j\in [1,\ell_{k+r}(0)]\cap K_\nu}\#\,(K_\nu\cap [j+1,j+\ell_k(1)])^N
 \\ &\le \frac{2^NJ}{\ell_{k+r}(0)}\max_{\e\in\Bbb Z_Q}\sum_{\nu\in\mathbb Z^d}V_k(0,\e)(\nu)^{N+1}
 \\ &\le\frac{2^N}{\ell_{k}(1)}\max_{\e\in\Bbb Z_Q}\sum_{\nu\in\mathbb Z^d}V_k(0,\e)(\nu)^{N+1}.\ \ \ \CheckedBox\text{\rm (4.2)}
 \end{align*}

\demo{Proof of the Visit Lemma}
 
 Let
 $$\widehat{V}_k(i,\e)(\th):=\sum_{\nu\in\mathbb Z^d}V_k(i,\e)(\nu)e^{2\pi\frak i \<\th,\nu\>}\ \ \ \  (\e\in\Bbb Z_Q,\ i=0,1,\ \th\in\Bbb T^d),$$
 then
 $$\widehat{V}_k(i,\e)(\th)=\ell_k(i)E(e^{2\pi\frak i \<\th, X^{(k)}(i,\e)\>}).$$
 \

 Using (4.1) in the sublemma and the Riesz-Fischer theorem, we see that

\begin{align*}\int_0^1\Psi_{\ell_k(0)}(x)dx\ & \ge\ \ 
\frac1{3\ell_k(0)}\min_{\e\in\Bbb Z_Q}\int_{\mathbb T^d}|\widehat{V}_k(1,\e)(\th)|^2d\th-\frac12\\ &=
\frac{\ell_k(1)^2}{3\ell_k(0)}\min_{\e\in\Bbb Z_Q}\int_{\mathbb T^d}|E(e^{2\pi\frak i \<\th, X^{(k)}(1,\e)\>})|^2d\th-\frac12.\end{align*}
This is (a).  To see (b),
  \begin{align*}\|\Psi_{\ell_k(1)}\|_\infty &\xleftarrow[N\to\infty]{}\ \(\int_0^1\Psi_{\ell_k(1)}(x)^Ndx\)^{\frac1N}\\ &\le
  \(\frac{2^N}{\ell_k(1)}\max_{\e\in\Bbb Z_Q}\sum_{\nu\in\mathbb Z^d}V_{k}(0,\e)(\nu)^{N+1}\)^{\frac1N}\ \ \text{by (4.2)}
  \\ &=
  \frac2{\ell_k(1)^{\frac1N}}\,\,\max_{\e\in\Bbb Z_Q}\(\sum_{\nu\in\mathbb Z^d}V_{k}(0,\e)(\nu)^{N+1}\)^{\frac1N}
  \\ &\le
 \frac2{\ell_k(1)^{\frac1N}}\,\,\,\,\max_{\e\in\Bbb Z_Q}\int_{\mathbb T^d}|\widehat{V}_{k}(0,\e)(\th)|^{1+\frac1N}d\th\ \ \text{by the Hausdorff-Young theorem}\\ &
  \xrightarrow[N\to\infty]{}\ 2\max_{\e\in\Bbb Z_Q}\int_{\mathbb T^d}|\widehat{V}_{k}(0,\e)(\th)|d\th\\ &=
  2\ell_k(0)\max_{\e\in\Bbb Z_Q}\int_{\mathbb T^d}|E(e^{2\pi\frak i \<\th, X^{(k)}(0,\e)\>})|d\th.
  \end{align*}
  
  This is (b).\ \ \CheckedBox
  
  \subsection*{Adaptedness}
   \ \ As in \cite{ABN2016},  the norm of a matrix $A\in M_{S\x S}$ is given by 
   $$\|A\|:=\sup\,\{\|Ax\|_\infty:\ x\in\Bbb R^S,\ \|x\|_\infty=1\}$$
   where $\|(x_s:\ s\in S)\|_\infty:=\sup_{s\in S}\,|x_s|$.

   We'll call  the  {\tt RAT} $F\in \text{\tt RV}(M_{S\x S}(\mathbb R)\x(\mathbb R^d)^S)$ {\it adapted} 
  if $\exists$ a discrete subgroup $\G=\G_{F}\le\mathbb R^d$ (called the {\it adaptivity group})  so that 
 $$\th\in\Bbb R^d\ \&\ \|\Pi_{F}(\th)\|=1\ \Lra\ \th\in\G.$$
 Equivalently, for some $r>0$,
 $$\|\G_F(\th)\|<1\ \forall\ \th\in B(0,r)\setminus\{0\}.$$
 Now, writing $F=(\mathcal L,W)\in\text{\tt RV}\,(S^S,(\Bbb R^d)^{S\x S})$, we have as $\th\to 0$
\begin{align*}\|\G_F(\th)\|&=\max_{s\in S}\sum_{t\in S}P(\mathcal L_s=t)|E(e^{2\pi\frak i \<\th, W_{s,t}\>})|\\ &
=1-\min_{s\in S}\sum_{t\in S}P(\mathcal L_s=t) \<\text{\tt Cov}\,(W_{s,t})\th,\th\>+o(\|\th\|)
\end{align*}
where for $V=(V_1,V_2,\dots,V_d)$ a $\Bbb R^d$-valued $L^2$ random variable, the {\it covariance matrix}
$\text{\tt Cov}\,(V)\in M_{d\x d}(\Bbb R)$ is defined by
$$\text{\tt Cov}\,(V)_{k,\ell}:=E((V_k-E(V_k))(V_\ell-E(V_\ell))).$$ 
\

A covariance matrix is {\it non-negative definite} in the sense that
$$\<\text{\tt Cov}\,(V)\th,\th\>=E\left((\sum_{k=1}^d\th_k(V_k-E(V_k)))^2\right)\ge 0\ \forall\ \th\in\Bbb R^d$$
and is called {\it positive definite} if it is invertible. Equivalently, for some $\e>0$ (the minimum eigenvalue modulus)
$$\<\text{\tt Cov}\,(V)\th,\th\>\ge \e\|\th\|^2\ \forall\ \th\in\Bbb R^d.$$
\

Thus, $F$ is adapted iff
                         
\begin{align*}\ \forall\ s\in S,\ &\exists\ t\in S\ \st\ P(\mathcal L_s=t)>0\ \&\\ & \text{\tt Cov}\,(W_{s,t})\ \text{ is strictly positive definite}. 
\end{align*}

It follows as in \cite{ABN2016} that if $F$ is adapted, then
$\forall\ \e>0\ \&\ M>0,\ \exists\ \d>0$ so that 
$$\|\G_{F}(\th)\|\le 1-\d\ \forall\ \ \ n\ge 1\ \ \text{and} \ \th\in B(0,M)\setminus B(\G,\e).$$

The following lemma gives a sequence version of adaptedness similar to that in \cite{ABN2016}.
   \proclaim{Adaptedness lemma}\ \ For $N,\ J,\ M\ge 1$ large, $\exists$
   a discrete subgroup $\G\le\Bbb R^d$ and $\e,\ b,\ c>0,\ r\in (0,1)$ so that
   \begin{align*}&\tag{i}\|\Pi_{\mathcal E_n}(\th+\g)\|\le 1-c\|\th\|^2\ \forall\ \g\in\G\cap B(0,M),\ \th\in B(0,r);
   \\ &\tag{ii}\|\Pi_{\mathcal E_n}(\th)\|\le 1-\e\ \forall\ \th\in B(0,M)\setminus\ B(\G,r);\\ &
   \tag{iii} \<\text{\tt Cov}\,(W_{s,t}(\mathcal E_n))\th,\th\>\ge\e\|\th\|^2\ \forall\ \th\in\Bbb R^d;
      \end{align*}
where
   $\mathcal E_n:=\mathcal F_{K+L(J+n)+1}^{K+L(J+n+1)}$ where $(\F_n\ :\ n\geq 1)$ is the independent {\tt RAT} sequence as on page \pageref{page: fancy F}.
   
   \endproclaim
 The proof is  in a series of steps, the first two of which are as in \cite{TM}.
 \Par1 If $\th\in\Bbb R^d,\ \xi\in\Bbb C\ \&\ v\in\Bbb C^S$ satisfy $\Pi_{\mathcal H}(\th)v=\xi v$, then
 $|\xi|\le 1$ with equality iff
  \begin{align*}\tag{a}v\in (\Bbb S^1)^S\ \&\  \ \
   (\Pi_{\mathcal H}(\th))_{s,t}=\xi v_s\overline{v}_t(\Pi_{\mathcal H}(0))_{s,t}\ \ \forall\ s,t\in S
  \end{align*}
  where $\Bbb S^1:=\{z\in \Bbb C\ :\ |z|=1\}$ is the multiplicative circle.
\demo{Proof} \ Write $\Pi:=\Pi_{\mathcal H}(\th)\ \&\ P:=\Pi_{\mathcal H}(0)$. Evidently (a) $\Lra\ \Pi v=\xi v$.
\

Now suppose that $\Pi v=\xi  v$ with $J\in S,\ |v_J|=\|v\|_\infty=1$, then
  \begin{align*}|\xi |&=|\xi  v_J|=
  \left|\sum_{t\in S}\Pi_{J,t}v_t\right|\le 
  \sum_{t\in S}|\Pi_{J,t}|\,|v_t|\\ &\le
  \sum_{t\in S}P_{J,t}\,|v_t|\le 1.   
  \end{align*}
If $|\xi |=1$, then $v\in (\Bbb S^1)^S$ and 
$$\left|\sum_{t\in S}\Pi_{s,t}v_t\right|=1\ \forall\ s\in S$$
and $\exists\ z\in (\Bbb S^1)^S$ so that
$$\Pi_{s,t}v_t=z_sP_{s,t}\ \forall\ s,t\in S.$$
Thus, for $s\in S$,
$$\xi  v_s=\sum_{t\in S}\Pi_{s,t}v_t=\sum_{t\in S}z_sP_{s,t}=z_s$$
which is (a).\ \Checkedbox

\

Next, for $\th\in\Bbb R^d$, let 
$$\mu(\th):=\max\,\{|\xi|:\ \xi\in\Bbb C\ \&\ \exists\ v\in\Bbb C^S,\ \Pi_{\mathcal H}(\th)=\xi v\}.$$
By \P1, $\mu(\th)\le 1$. Set $\G:=\{\g\in\Bbb R^d:\ \mu(\g)=1\}$.
\Par2 $\G$ is a discrete subgroup of $\Bbb R^d$ and
\begin{align*}
 \tag{b}\mu(\g+\th)=|\l(\th)|\ \forall\ \g\in\G,\ \|\th\|_\infty\le r_{\mathcal H}.
\end{align*}
\demo{Proof}\ \ Since $P(\mathcal L_s=t)>0\ \forall\ s,t\in S$,
$$\G=\{\g\in\Bbb R^d:\ \exists\ x\in\Bbb R,\ \<\g, W_{s,t}\>+x\in 2\pi\Bbb Z\ \text{a.s.}\ \forall\ s,t\in S\}$$
which is evidently a subgroup of $\Bbb R^d$.
\

Now suppose that $\g\in\G$ with $\Pi_\mathcal H(\g)v=\xi v$ where $|v_s|=|\xi|=1\ \forall\ s\in S$.
By \P1,
$$(\Pi_{\mathcal H}(\th))_{s,t}=\xi v_s\overline{v}_t(\Pi_{\mathcal H}(0))_{s,t}\ \ \forall\ s,t\in S.$$
Equivalently $\forall\ s,t\in S,$
$$E(e^{2\pi\frak i\<\g,W_{s,t}\>})=\xi v_s\overline{v}_t\ \Lra\  E(e^{2\pi\frak i\<(\g+\th),W_{s,t}\>})=\xi v_s\overline{v}_tE(e^{2\pi\frak i \<\th, W_{s,t}\>})
\ \forall\ \th\in\Bbb R^d$$
whence
$$(\Pi_{\mathcal H}(\th+\g))_{s,t}=\xi v_s\overline{v}_t(\Pi_{\mathcal H}(\th))_{s,t}\ \ \forall\ s,t\in S,\ \th\in\Bbb R^d.$$
 Statement (b) follows from this, whence \P2 via the Taylor expansion of $\l$.\ \Checkedbox
 
  \Par3 For $N\ge 1$ sufficiently large, $\mathcal H_1^N$ is an adapted {\tt RAT} with adaptivity group $\G$.
     \demo{Proof} \ Fix $0<q<1$ and $N\ge 1$ so that
   
   \begin{align*}&\|Q(\th)^N\|<q\ \text{for}\ \th\in B(0,r_\mathcal H)\ \&\ \text{hence for}\ \th\in B(\G,r_\mathcal H); \ \&
   \\ & \mu(\th)^N<q\ \forall\ \th\in\Bbb R^d\setminus B(\G,r_\mathcal H). \end{align*}
   It follows that $\mathcal H_1^N$ is  adapted with adaptivity group $\G$.\ \ \Checkedbox
 
 To complete the proof of the lemma, fix  $J\ge 1$ and let 
 $$(\mathcal E_n:=\mathcal F_{J+Nn+1}^{J+N(n+1)}:\ n\ge 1).$$
  Statements (i) and (ii) follow because for each $M>0$,
 $$\sup_{\|\th\|\le M}\|\Pi_{\mathcal E_n}(\th)-\Pi_{\mathcal H_1^N}(\th)\|\xrightarrow[n\to\infty]{}\ \ 0$$
 and statement (iii) follows from
 $$W_{s,t}(\mathcal E_n)\xrightarrow[n\to\infty]{\text{\tt\tiny RV}\,(\Bbb R^d)}\ W_{s,t}(\mathcal H_1^N)\ \forall\ s,t\in S.\ \CheckedBox$$

     \

    To establish  theorem 3, we use the following.
    \proclaim{Weak, rough local limit theorem}\ \ \ \ 
 
For each $s\in S$ and $1\le p\le 2$,
\begin{align*}\tag{\tt WRLLT}
 \int_{\mathbb T^d}|E(e^{2\pi\frak i \<\th,X^{(J+Ln)}_s\>})|^pd\th\ \text{\Large $\asymp$}\ \frac1{n^\frac{d}2}.
\end{align*}

\endproclaim
The proof is a multidimensional version of the proof of the {\tt WRLLT} in \cite{ABN2016}.
    \demo{Proof of $\gg$}

    \

By theorem 6.1 in  \cite{ABN2016}, for each $1\le k\le d$, we have
    $$\sum_{\nu=1}^n\min_{s\in S}E([(b_s(\mathcal E_\nu))_k]^2)\le E([(X_s^{(J+Ln)})_k]^2)\le\sum_{\nu=1}^n\max_{s\in S}E([(b_s(\mathcal E_\nu))_k]^2).$$
Thus, by the Adaptedness lemma,$\exists\ G>0$ so that
$$E(\|X_s^{(J+Ln)}\|^2)\le G n.$$

Next, fix $M=2\sqrt{G}$, then by Chebyshev's inequality,

$$P([\|X_s^{(J+Ln)}\|\le M\sqrt n])\ge\frac34.$$
Now fix $\D>0$ so that 
$$|1-e^{2\pi\frak ix}|<\frac14\ \forall\ |x|<\D.$$
We have
\begin{align*}n^{\frac{d}2}\int_{[-\pi,\pi]^d}|E(e^{2\pi\frak i\<\th,X^{(J+Ln)}_s\>})|^2d\th &=\int_{[-\pi\sqrt n,\pi\sqrt n]^d}|E(\exp[2\pi\frak i \<\th,\frac{X^{(J+Ln)}_s}{\sqrt n}\>])|^2d\th\\ &\ge
\int_{[-\frac{\D}M,\frac{\D}M]^d}|E(\exp[2\pi\frak i \<\th,\frac{X_s^{(J+Ln)}}{\sqrt n}\>])|^2d\th.
 \end{align*}
For $\|\th\|<\frac{\D}M$, we have
\begin{align*}
 |E(\exp[2\pi\frak i \<\th,&\frac{X^{(J+Ln)}_s}{\sqrt n}\>])|\\ &\ge |E(\exp[2\pi\frak i \<\th,\frac{X^{(J+Ln)}_s}{\sqrt n}\>])1_{[\|X_s^{(J+Ln)}\|<M\sqrt n]})|-P([\|X_s^{(J+Ln)}\|\ge M\sqrt n]\\ &\ge 
 \frac34\cdot\frac34-\frac14\\ &=\frac{5}{16}\end{align*}
 whence
 \begin{align*}n^{\frac{d}2}\int_{[-\pi,\pi]^d}|E(e^{2\pi\frak i\<\th,X^{(J+Ln)}_s\>})|^2d\th &\ge
\int_{[-\frac{\D}M,\frac{\D}M]^d}|E(\exp[2\pi\frak i \<\th,\frac{X^{(J+Ln)}_s}{\sqrt n}\>])|^2d\th\\ &\ge \(\frac{2\D}M\)^d\cdot\frac{25}{256}.\ \ \ \CheckedBox\ \ \gg
 \end{align*}
 \demo{Proof of $\ll$}
 \

  We have
\begin{align*}
 |E(e^{2\pi\frak i\<\th,X_s^{(J+Ln)}\>})|&=
 |(\Pi_{\mathcal E_n}(\th)\Pi_{\mathcal E_{n-1}}(\th)\cdots \Pi_{\mathcal E_1}(\th)\widehat{\Xi}_J(\th))_s|\\ &
 \le \prod_{k=1}^n\|\Pi_{\mathcal E_k}(\th)\|.
\end{align*}
Fix $M>0$ so that $[-\pi,\pi]^d\subset B(0,M)$. By the Adaptedness lemma, for $n\ge 1,\ \g\in\G$ we have
$$\|\Pi_{\mathcal E_n}(\g+\th)\|\le 1-c\|\th\|^2\ \ \forall\ |\th|<r$$ and
$$\|\Pi_{\mathcal E_n}(\th)\|\le 1-\e\ \ \forall\ \th\in B(0,M)\setminus B(\G,r).$$
\begin{align*}\int_{[-\pi,\pi]^d}&|E(e^{2\pi\frak i\<\th,X_s^{(J+Ln)}\>})|d\th \\ &\le 
(\int_{B(0,M)\cap B(\G,r)}+\int_{B(0,M)\setminus B(\G,r)}) \prod_{k=1}^n\|\Pi_{\mathcal E_k}(\th)\|d\th\\ &\le
\sum_{\g\in B(0,M)\cap \G}\int_{B(0,r)}\prod_{k=1}^n\|\Pi_{\mathcal E_k}(\g+\th)\|d\th+\int_{B(0,M)\setminus B(\G,r)}\prod_{k=1}^n\|\Pi_{\mathcal E_k}(\th)\|d\th\\ &\le
\#\,(B(0,M)\cap \G)\int_{B(0,r)}(1-c\|\th\|^2)^nd\th+O((1-\e)^n)\\ &
\ll\frac1{n^\frac{d}2}.\ \ \ \CheckedBox\ \ll\ \&\ \text{\tt WRLLT} \end{align*}

     \demo{Proof of theorem 3}
 \
 
Set $\nu_k(i):=\ell_{K+Lk}(i)\ \ (i=0,1)$. The Visit lemma  and the {\tt WRLLT} show that
  \begin{align*}\Psi_{\nu_k(i)}(x)&\ll \nu_k(0)\int_{\Bbb T^d}|E(e^{2\pi\frak i\<\th,X_{(0,i)}^{(J+Lk)}\>})|d\th\\ &\asymp
  \frac{\nu_k(0)}{k^{\frac{d}2}}\ll
  \nu_k(0)\int_{\Bbb T^d}|E(e^{2\pi\frak i\<\th,X_s^{(J+Lk)}\>})|^2d\th\\ &\asymp\int_{\Bbb T^d}\Psi_{\nu_k(i)}(x)dx\\ &\ll
  \frac{\nu_k(0)}{k^{\frac{d}2}}   
  \end{align*}
  \

  Next, $\exists\ \Lambda>1$ so that
$\nu_k(i)\propto \Lambda^k\ \ (i=0,1)$
whence
for \par $\nu_k(0)\le n\le \nu_{k+1}(0)$,
\begin{align*}
 \int_{\mathbb T^d}\Psi_{n}(x)dx&\ge \int_{\mathbb T^d}\Psi_{\nu_k(0)}(x)dx\\& \gg\ 
\frac{{\nu_k}(1)}{k^{\frac{d}2}}\gg\frac{\nu_{k+1}(0)}{k^{\frac{d}2}}\\ &\gg
\frac{n}{(\log n)^{\frac{d}2}}
\end{align*}

and for ${\nu_k}(1)\le n\le \nu_{k+1}(1)$,
\begin{align*}
 \|\Psi_n(x)\|_{L^\infty(\Bbb T^d)}\ll   \frac{{\nu_{k+1}}(1)}{k^{\frac{d}2}}\ll\frac{{\nu_k}(0)}{k^{\frac{d}2}}\ll \frac{n}{(\log n)^{\frac{d}2}}.\ \ \ \CheckedBox
\end{align*}


\begin{bibdiv}
\begin{biblist}

\bib{ABN2016}{article}{
      author={{Aaronson}, Jon.},
      author={{Bromberg}, Michael},
      author={{Nakada}, Hitoshi},
       title={{Discrepancy Skew Products and Affine Random Walks}},
        date={2016-03},
     journal={ArXiv e-prints},
      eprint={1603.07233},
}

\bib{AK}{article}{
      author={Aaronson, Jon.},
      author={Keane, Michael},
       title={The visits to zero of some deterministic random walks},
        date={1982},
        ISSN={0024-6115},
     journal={Proc. London Math. Soc. (3)},
      volume={44},
      number={3},
       pages={535\ndash 553},
         url={http://dx.doi.org/10.1112/plms/s3-44.3.535},
      review={\MR{656248 (83j:28016)}},
}

\bib{Beck}{book}{
      author={Beck, J{\'o}zsef},
       title={Probabilistic {D}iophantine approximation},
      series={Springer Monographs in Mathematics},
   publisher={Springer, Cham},
        date={2014},
        ISBN={978-3-319-10740-0; 978-3-319-10741-7},
         url={http://dx.doi.org/10.1007/978-3-319-10741-7},
        note={Randomness in lattice point counting},
      review={\MR{3308897}},
}

\bib{BU}{article}{
      author={{Bromberg}, Michael},
      author={{Ulcigrai}, Corinna},
       title={{A temporal Central Limit Theorem for real-valued cocycles over
  rotations}},
        date={2017-05},
     journal={ArXiv e-prints},
      eprint={1705.06484},
}

\bib{Conze76}{incollection}{
      author={Conze, Jean-Pierre},
       title={Equir\'epartition et ergodicit\'e de transformations
  cylindriques},
        date={1976},
   booktitle={S\'eminaire de {P}robabilit\'es, {I} ({U}niv. {R}ennes, {R}ennes,
  1976), {E}xp. {N}o. 2},
   publisher={D\'ept. Math. Informat., Univ. Rennes, Rennes},
       pages={21},
      review={\MR{0584019}},
}

\bib{CP}{article}{
      author={Conze, Jean-Pierre},
      author={Pi{\c{e}}kniewska, Agata},
       title={On multiple ergodicity of affine cocycles over irrational
  rotations},
        date={2014},
        ISSN={0021-2172},
     journal={Israel J. Math.},
      volume={201},
      number={2},
       pages={543\ndash 584},
         url={http://dx.doi.org/10.1007/s11856-014-0033-3},
      review={\MR{3265295}},
}

\bib{DS}{article}{
      author={Dolgopyat, Dmitry},
      author={Sarig, Omri},
       title={Temporal distributional limit theorems for dynamical systems},
        date={2017},
        ISSN={0022-4715},
     journal={J. Stat. Phys.},
      volume={166},
      number={3-4},
       pages={680\ndash 713},
         url={http://dx.doi.org/10.1007/s10955-016-1689-3},
      review={\MR{3607586}},
}

\bib{HW}{book}{
      author={Hardy, G.~H.},
      author={Wright, E.~M.},
       title={An introduction to the theory of numbers},
   publisher={Oxford, at the Clarendon Press},
        date={1954},
        note={3rd ed},
      review={\MR{0067125}},
}

\bib{HH}{book}{
      author={Hennion, Hubert},
      author={Herv{\'e}, Lo{\"{\i}}c},
       title={Limit theorems for {M}arkov chains and stochastic properties of
  dynamical systems by quasi-compactness},
      series={Lecture Notes in Mathematics},
   publisher={Springer-Verlag, Berlin},
        date={2001},
      volume={1766},
        ISBN={3-540-42415-6},
         url={http://dx.doi.org/10.1007/b87874},
      review={\MR{1862393 (2002h:60146)}},
}

\bib{Herman}{article}{
      author={Herman, Michael-Robert},
       title={Sur la conjugaison diff\'erentiable des diff\'eomorphismes du
  cercle \`a des rotations},
        date={1979},
        ISSN={0073-8301},
     journal={Inst. Hautes \'Etudes Sci. Publ. Math.},
      number={49},
       pages={5\ndash 233},
         url={http://www.numdam.org/item?id=PMIHES_1979__49__5_0},
      review={\MR{538680}},
}

\bib{Katznelson}{article}{
      author={Katznelson, Y.},
       title={Sigma-finite invariant measures for smooth mappings of the
  circle},
        date={1977},
        ISSN={0021-7670},
     journal={J. Analyse Math.},
      volume={31},
       pages={1\ndash 18},
      review={\MR{0486415}},
}

\bib{Keane1970}{incollection}{
      author={Keane, Michel},
       title={Irrational rotations and quasi-ergodic measures},
        date={1970},
   booktitle={Publications des {S}\'eminaires de {M}ath\'ematiques ({U}niv.
  {R}ennes, {R}ennes, ann\'ee 1970--1971), {F}asc. 1: {P}robabilit\'es},
   publisher={D\'ep. Math. et Informat., Univ. Rennes, Rennes},
       pages={17\ndash 26},
      review={\MR{0369659}},
}

\bib{Khintchine}{book}{
      author={Khintchine, A.~Ya.},
       title={Continued fractions},
      series={Translated by Peter Wynn},
   publisher={P. Noordhoff, Ltd., Groningen},
        date={1963},
      review={\MR{0161834}},
}

\bib{KN}{article}{
      author={Kraaikamp, Cor},
      author={Nakada, Hitoshi},
       title={On normal numbers for continued fractions},
        date={2000},
        ISSN={0143-3857},
     journal={Ergodic Theory Dynam. Systems},
      volume={20},
      number={5},
       pages={1405\ndash 1421},
         url={http://dx.doi.org/10.1017/S0143385700000766},
      review={\MR{1786721 (2001i:11101)}},
}

\bib{Kronecker}{article}{
      author={Kronecker, L.},
       title={Zwei {S}\"atze \"uber {G}leichungen mit ganzzahligen
  {C}oefficienten},
        date={1857},
        ISSN={0075-4102},
     journal={J. Reine Angew. Math.},
      volume={53},
       pages={173\ndash 175},
         url={http://dx.doi.org/10.1515/crll.1857.53.173},
      review={\MR{1578994}},
}

\bib{Oren}{article}{
      author={Oren, Ishai},
       title={Ergodicity of cylinder flows arising from irregularities of
  distribution},
        date={1983},
        ISSN={0021-2172},
     journal={Israel J. Math.},
      volume={44},
      number={2},
       pages={127\ndash 138},
         url={http://dx.doi.org/10.1007/BF02760616},
      review={\MR{693356}},
}

\bib{rigveda}{book}{
      author={Schmidt, Klaus},
       title={Cocycles on ergodic transformation groups},
   publisher={Macmillan Company of India, Ltd., Delhi},
        date={1977},
        note={Macmillan Lectures in Mathematics, Vol. 1},
      review={\MR{0578731}},
}

\bib{TM}{incollection}{
      author={Taussky, Olga},
      author={Marcus, Marvin},
       title={Eigenvalues of finite matrices: {S}ome topics concerning bounds
  for eigenvalues of finite matrices},
        date={1962},
   booktitle={Survey of numerical analysis},
   publisher={McGrawHill, New York},
       pages={279\ndash 297},
      review={\MR{0132071}},
}

\end{biblist}
\end{bibdiv}
\end{document}